\documentclass{amsart}

\usepackage[utf8]{inputenc}
\usepackage[T1]{fontenc}

\usepackage{geometry}
\usepackage{ragged2e}

\usepackage[pagebackref, hidelinks]{hyperref}
\hypersetup{
colorlinks=true, 
citecolor = blue, 
linkcolor = blue, 
urlcolor = blue, 
anchorcolor = blue
}

\usepackage{amsmath, amsthm, amsfonts}
\usepackage{amscd, amssymb, latexsym}

\usepackage[all]{xy}
\usepackage{mathtools} 
\usepackage{caption}

\usepackage{paratext/commands_LiNuMoK}

\usepackage{enumitem}
\usepackage{comment}
\usepackage{color}

\usepackage{diagbox} 
\usepackage{tikz, pgf, pgfplots}
\pgfplotsset{compat=newest}

\usetikzlibrary[patterns]
\usetikzlibrary{arrows, arrows.meta}
\usetikzlibrary{calc} 
\def\centerarc[#1](#2)(#3:#4:#5)%
    {\draw[#1] ($(#2)+({#5*cos(#3)},{#5*sin(#3)})$) arc (#3:#4:#5);}

\usetikzlibrary{intersections,through,backgrounds} 
\usetikzlibrary{shapes.geometric} 

\definecolor{grey}{rgb}{0.5,0.5,0.5}
\definecolor{forestgreen}{rgb}{0.133,0.545,0.133}
\definecolor{marron}{rgb}{0.6,0.2,0.}
\definecolor{violet}{rgb}{0.5,0.,0.5}
\definecolor{lightpurple}{rgb}{.6,.2,.6}

\title{Linking numbers of modular knots}
\author{Christopher-Lloyd Simon}
\date{\today}

\begin{document}

\maketitle

\begin{abstract}
    The modular group $\PSL_2(\Z)$ acts on the hyperbolic plane $\H\P$ with quotient the modular surface $\M$, whose unit tangent bundle $\U$ is a $3$-manifold homeomorphic to the complement of the trefoil knot in the sphere $\S^3$.
    The hyperbolic conjugacy classes of $\PSL_2(\Z)$ correspond to the closed oriented geodesics in $\M$. Those lift to the periodic orbits for the geodesic flow in $\U$, which define the modular knots. 
    
    The linking numbers between modular knots and the trefoil are well understood.
    Indeed, \'E. Ghys showed in 2006 that they are given by the Rademacher invariant of the corresponding conjugacy classes. The Rademacher function $\Rad \colon \PSL_2(\Z)\to \Z$ is a homogeneous quasi-morphism which he recognised with J. Barge in 1992 as half the primitive of the bounded Euler class.
    
    We are concerned with the linking numbers between modular knots and derive several formulae with arithmetical, combinatorial, topological and group theoretical flavours.
    In particular we associate to a pair of modular knots a function defined on the character variety of $\PSL_2(\Z)$, whose limit at the boundary point recovers their linking number.
    Moreover, we show that the linking number with a modular knot minus that with its inverse yields a homogeneous quasi-morphism on the modular group, and how to extract a free basis out of these. 
\end{abstract}

\renewcommand{\contentsname}{Plan of the paper}
\setcounter{tocdepth}{1}
\tableofcontents

\subsection*{Acknowlegements}
This paper contains the main results obtained in the second part of my thesis.
I would thus like to thank my thesis advisors Etienne Ghys and Patrick Popescu-Pampu for their guidance and encouragement; as well as Francis Bonahon, Louis Funar, Jean-Pierre Otal and Anne Pichon who refereed and carefully read my work.
I am also grateful to Pierre Dehornoy for sharing his knowledge of modular knots, as well as to the members of the arithmetics and dynamics teams at PSU for inviting me to present my works in their seminars.
I owe Marie Dossin for helping me with the figures in tikz.
Finally, I thank the referee for their relevant suggestions and questions.

\section{Introduction}

\subsection*{Context and motivation}

The modular group $\PSL_2(\Z)$ acts properly discontinuously on the hyperbolic plane $\H\P$ with quotient the modular orbifold $\M$, a hyperbolic surface with conical singularities $i$ and $j$ of order $2$ and $3$, and a cusp $\infty$.
The free homotopy classes of loops in $\M$ correspond to the conjugacy classes of its fundamental group $\pi_1(\M)=\PSL_2(\Z)$.
In particular the hyperbolic conjugacy classes of $\PSL_2(\Z)$ correspond to the closed oriented geodesics in $\M$, called \emph{modular geodesics}.
For hyperbolic $A\in \PSL_2(\Z)$ the modular geodesic $\alpha$ has length $\lambda_A$ equal to the logarithm of the ratio between its eigenvalues $\epsilon_A^{\pm 1}$, in formula:
\begin{equation*}
    \disc(A) 
    = \left(\epsilon_A-\epsilon_A^{-1}\right)^2 
    = (\Tr A)^2-4
    = \left(2\sinh \tfrac{1}{2}\lambda_A\right)^2
\end{equation*}
We denote by $I(A,B)$ the geometric intersection number between the associated modular geodesics.

The unit tangent bundle $\U=\PSL_2(\Z)\backslash \PSL_2(\R)$ of $\M=\PSL_2(\Z)\backslash \H\P$ is a $3$-manifold, and the closed oriented geodesics in $\M$ lift to the periodic orbits for the geodesic flow in $\U$.
Hence the primitive hyperbolic conjugacy classes in $\PSL_2(\Z)$ correspond to the so called \emph{modular knots} in $\U$, which form the components of the \emph{master modular link}.
The structure of the Seifert fibration $\U \to \M$ reveals that $\U$ is homeomorphic to the complement of a trefoil knot in the sphere. In particular, one may speak of the linking numbers between modular knots and the trefoil and between one another.

Let us recall a combinatorial parametrization of the infinite order conjugacy classes in $\PSL_2(\Z)$.
The Euclidean algorithm shows that the group $\SL_2(\Z)$ is generated by the transvections
\begin{equation*}
    L=
    \begin{psmallmatrix}
    1 & 0 \\ 1 & 1
    \end{psmallmatrix} 
    \qquad \mathrm{and} \qquad
    R=
    \begin{psmallmatrix}
    1 & 1 \\ 0 & 1
    \end{psmallmatrix}
\end{equation*}
and more precisely that its submonoid $\SL_2(\N)$ of matrices with non-negative entries is freely generated by $\{L,R\}$.
This submonoid can be identified with its image $\PSL_2(\N)\subset \PSL_2(\Z)$.
In $\PSL_2(\Z)$, the conjugacy class of an infinite order element intersects $\PSL_2(\N)$ along all cyclic permutations of a non-empty $\{L,R\}$-word.
The conjugacy class is primitive if and only if the cyclic word is primitive, and it is hyperbolic when the cyclic word contains both letters $L$ and $R$.

One may try to relate the geometry and topology of the master modular link with the arithmetics and combinatorics of conjugacy classes in the modular group.
Our previous work \cite{CLS_Conj-PSL2K_2022} relates the geometry of modular geodesics (angles of intersection and lengths of ortho-geodesics) in terms of the arithmetics of conjugacy classes in the modular group (discriminants, cross-ratios between fixed points on the projective line, and their Hilbert symbols).
The main results in this paper will relate the linking numbers of modular knots to the combinatorics of the corresponding cyclic words.

The most immediate measures for the complexity of a binary word are given by the sum and difference between the numbers of letters of each sort.
For an infinite order $A\in \PSL_2(\Z)$ we denote by $\len([A])=\#R+\#L$ the combinatorial length and call $\Rad([A])=\#R-\# L$ the \emph{Rademacher number} of its conjugacy class.
In his paper \cite{Atiyah_log(eta-Dedekind)_1987} on the Logarithm of the Dedekind eta function, M. Atiyah identified the Rademacher function with no less than six other important functions appearing in diverse areas of mathematics, showing how omnipresent it is.

The function $\Rad\colon \PSL_2(\Z) \to \Z$ is a quasi-morphism, meaning that it has a bounded derivative
\begin{equation*}
d\Rad\colon \PSL_2(\Z)\times \PSL_2(\Z)\to \Z 
\qquad 
d\Rad(A,B)=\Rad(B)-\Rad(AB)+\Rad(A)
\end{equation*}
and is homogeneous, meaning that $\Rad(A^n)=n\Rad(A)$ for infinite order $A\in \PSL_2(\Z)$ and $n\in \Z$.
This enabled \'E. Ghys and J. Barge to recognise it in \cite{BargeGhys_cocycle-euler-maslov_1992} as half the primitive of the bounded Euler class in $H^2_b(\PSL_2(\Z);\R)$ and explain its ubiquity.
In \cite{Ghys_knots-dynamics_2006}, \'E. Ghys showed that the linking number of a modular knot with the trefoil equals its Rademacher invariant, and concluded by asking for \emph{arithmetical and combinatorial interpretations of the linking pairing between modular knots}.

In this work, we will derive several formulae for those linking numbers, providing bridges between the arithmetics and geometry, the combinatorics and algebra, or the dynamics and topology of the modular group.

\subsection*{The arithmetic and geometry of the cosines}

The journey from arithmetics began with our previous work \cite{CLS_Conj-PSL2K_2022}. In particular for a field $\Field\supset \Q$, we described when two hyperbolic elements $A,B\in \PSL_2(\Z)$ of the same discriminant $\Delta$ are conjugate in $\PSL_2(\Field)$.
The obstruction is measured in terms of any intersection angle $\theta$ between their modular geodesics $[\alpha],[\beta]$ by the class of $\left(\cos \tfrac{\theta}{2}\right)^2 \in \Field^\times$ modulo the group of norms over $\Field$ of the extension $ \Field[\sqrt{\Delta}]$.
When this obstruction vanishes, the elements of $\PSL_2(\Field)$ which conjugate $A$ to $B$ are paremetrized by the points $(X,Y)\in \Field$ of the generalised Pell-Fermat conic with equation $X^2-\Delta Y^2 =  \left(\cos \tfrac{\theta}{2}\right)^2$.
Hence the geometric quantities $\left(\cos \tfrac{\theta}{2}\right)^2=\tfrac{1+\cos(\theta)}{2}$ given for some representatives $A,B\in \SL_2(\R)$ whose axes intersect in $\H\P$ by
\begin{equation*}
    \cos(\theta) = \sign(\Tr(A)\Tr(B))\tfrac{\Tr(AB)-\Tr(AB^{-1})}{\sqrt{\disc A\disc B}}
\end{equation*}
have an arithmetic meaning, and they will reappear under various forms in the sequel.

\subsection*{Linking functions on the character variety}

Let us introduce, for any pair of modular geodesics $[\alpha],[\beta]\subset \M$, the following summations over their oriented intersection angles $\theta \in \,]-\pi,\pi[$:
\begin{equation*}
    \Link_1(A,B)
    = \tfrac{1}{2} \sum \left(\cos \tfrac{\theta}{2}\right)^2 
    \qquad \mathrm{and} \qquad
    \Cos_1(A,B)
    = \tfrac{1}{2} \sum \left(\cos \theta\right)
\end{equation*}
and study their variations as we deform the metric on $\M$ by opening the cusp.

The complete hyperbolic metrics on the orbifold $\M$ correspond to the faithful and discrete representations $\rho\colon \PSL_2(\Z) \to \PSL_2(\R)$ up to conjugacy.
They form a $1$-dimensional real algebraic set parametrized by $q\in \R^*$, deforming $L$ and $R$ according to:
\begin{equation*}
    L_q =
    \begin{pmatrix}
    q & 0 \\
    1 & q^{-1}
    \end{pmatrix}
    \qquad \mathrm{and} \qquad
    R_q = 
    \begin{pmatrix}
    q & 1 \\
    0 & q^{-1}
    \end{pmatrix}.
\end{equation*}
The primitive hyperbolic conjugacy classes of $\PSL_2(\Z)$ still index the hyperbolic geodesics in the quotient $\M_q=\rho_q(\PSL_2(\Z))\backslash\H\P$ which do not surround the cusp. The intersection points between $q$-modular geodesics $[\alpha_q],[\beta_q]\subset \M_q$ persist by deformation, so we may define the analogous sums $\Link_q(A,B)$ and $\Cos_q(A,B)$ over their intersection angles $\theta_q \in \,]-\pi,\pi[$ of the $\tfrac{1}{2}\left(\cos \tfrac{1}{2}\theta_q\right)^{2}$ and $\tfrac{1}{2}\left(\cos \theta_q\right)$.

As $q\to \infty$, the hyperbolic orbifold $\M_q$ has a convex core which retracts onto a thin neighbourhood of the long geodesic arc connecting its conical singularities, whose preimage in the universal cover $\H\P$ is a trivalent tree. In the limit we recover the action of $\PSL_2(\Z)$ on its Bruhat-Tits building, the infinite planar trivalent tree $\Tree$, and by studying its combinatorics we shall prove the following.

\begin{Theorem}[Linking and intersection from boundary evaluations]
For primitive hyperbolic $A,B\in \PSL_2(\Z)$, the limits of the functions $\Link_q(A,B)$ and $\Cos_q(A,B)$ at the boundary point of the $\PSL_2(\R)$-character variety of $\PSL_2(\Z)$ recover their linking and intersection numbers:
\begin{align*}
    &\Link_q(A,B) \xrightarrow[q\to \infty]{} \lk(A,B)
    \\
    &\Cos_q(A,B) \xrightarrow[q\to \infty]{} 
    \lk(A,B)-\lk(A^{-1},B)
    =2\lk(A,B)-\tfrac{1}{2}I(A, B)
\end{align*}
\end{Theorem}

Hence the functions $\Link_q$ and $\Cos_q$ interpolate between the geometry at $q=1$ of the arithmetic group $\PSL_2(\Z) \subset \PSL_2(\R)$ and the topology at $q=\infty$ of the combinatorial action $\PSL_2(\Z) \to \Aut(\Tree)$.

\subsection*{Linking functions and Alexander polynomials}

The graphs of $q\in \C \mapsto \Link_q(A,B) \in \C$ for various pairs $A,B$ suggest that their zeros tend to accumulate on the unit circle. This reminds us of the various results and conjectures concerning the roots of Alexander polynomials, so we propose a possible thread to follow in this direction. 

A primitive hyperbolic conjugacy class in $\PSL_2(\Z)=\pi_1(\M)$ corresponds to a primitive modular geodesic in $\alpha\subset \M$. It lifts to a modular knot in $\Vec{\alpha} \subset \U$ which in turns yields a conjugacy class in $\BB_3=\pi_1(\U)$.
A conjugacy class in the braid group on three strands defines, by taking its closure, a link $\sigma_A$ in a solid torus.
In \cite[Proposition 5.16]{CLS_phdthesis_2022} we relate the Alexander polynomial $\Delta(\sigma_A)\in \Z[t^{\pm 1}]$ of this link $\sigma_A$ to the Fricke polynomial $\Tr A_q \in \Z[q^{\pm 1}]$ of the modular geodesic $[\alpha]$.
\begin{Proposition}[Alexander, Fricke and Rademacher]
For a primitive hyperbolic $A\in \SL_2(\N)$, the Alexander polynomial of the link $\sigma_A$ is given in terms of $q=\sqrt{-t}$ by: \[\Delta(\sigma_A)=\tfrac{q^{\Rad(A)}-\Tr(A_q)+q^{-\Rad(A)}}{(q-q^{-1})^2}\]
\end{Proposition}
Now recall that $\Cos_q(A,B)=\Link_q(A,B)-\Link_q(A,B^{-1})$ can be expressed as finite sum of terms:
\begin{equation*}
    \cos(A_q,B_q) = \sign(\Tr(A_q)\Tr(B_q))\tfrac{\Tr(A_qB_q)-\Tr(A_qB_q^{-1})}{\sqrt{\disc A_q\disc B_q}} 
    \qquad \mathrm{where} \qquad 
    \disc(C_q)= (\Tr C_q)^2-4
\end{equation*}
This is how one may compare the concentration property for the zeros of $\Link_q$ around the unit circle with those of Alexander polynomials, but we will not pursue this direction any further.

\subsection*{Linking numbers and homogeneous quasi-morphisms}

The limiting values $\Cos_q(A,B)\to \lk(A,B)-\lk(A^{-1},B)$ will now provide a bridge from the representation theory to the bounded cohomology of $\PSL_2(\Z)$.
For every group $\Pi$, the real vector space $PX(\Pi)$ of homogeneous quasi-morphisms is a Banach space for the norm $\lVert df \rVert_\infty$, as was shown in \cite{MatsuMorita_Hb(Homeo)_1985, Ivanov_H2b(G)-Banach_1988}.

\begin{Theorem}
For every hyperbolic $A\in \PSL_2(\Z)$, the function $\Cos_A\colon B\mapsto \lk(A,B)-\lk(A^{-1}, B)$ is a homogeneous quasi-morphism $\PSL_2(\Z)\to \Z$.
\end{Theorem}

Let $\mathcal{P}$ denote the set of primitive infinite order conjugacy classes in $\PSL_2(\Z)$, and $\mathcal{P}_0$ the subset of those which are stable under inversion. Choose a partition $\mathcal{P}\setminus \mathcal{P}_0=\mathcal{P}_-\sqcup \mathcal{P}_+$ in two subsets in bijection by the inversion.
We may choose $R\in \mathcal{P}_+$, and denote $\Cos_R:=\Rad$ by convention.

\begin{Theorem}
The collection of $\Cos_A\in PX(\PSL_2(\Z))$ for $A\in \mathcal{P}_+$ is linearly independent and every element $f\in PX(\PSL_2(\Z))$ can be written as $f=\sum_{A\in \mathcal{P}_+} c_A(f) .\Cos_A$ for unique $c_A(f) \in \R$.
\end{Theorem}

To prove the non-triviality and linear independence of the $\Cos_A$ for $A\in \mathcal{P}_+$, we were led to show the non-degeneracy of the linking form, which is interesting in its own right.

\begin{Theorem}
If hyperbolic $A,B\in \PSL_2(\Z)$ are link equivalent, namely $\lk(A,X)=\lk(B,X)$ for all hyperbolic $X\in \PSL_2(\Z)$, 
then they are conjugate.
\end{Theorem}

The results in this section may be compared to the classical representation theory of compact groups, in which the characters of irreducible representations provide an orthonormal basis for the class functions.
Indeed, we have found a family of cosign functions $\Cos_A$ whose periods correspond to the primitive conjugacy classes of $\PSL_2(\Z)$, and they form a basis for the space of quasi-characters $PX(\PSL_2(\Z))$, which is in some sense orthogonal with respect to the linking form (but we refer to the proofs in section \ref{sec:quasi-morphism} for a better explanation of this orthogonality).

\section{The group \texorpdfstring{$\PSL_2(\R)$}{PSL(2;R)}: discriminant and cross-ratio}
\label{sec:disc-bir}

The group $\PSL_2$ acts by conjugacy over itself.
In this section, we recall from \cite{CLS_Conj-PSL2K_2022} the main invariants which enable to describe the orbits of single elements and of pairs of elements.
Those are the discriminant $\disc(A)$ and the cross-ratio $\bir(A,B)$.

\subsection{Over a field of characteristic $\ne 2$}
\label{subsec:disc-bir_K}

Let $\Field$ be a field of characteristic different from $2$. 
The automorphism group $\PGL_2(\Field)$ of the projective line $\Field\P^1$ acts freely transitively on triples of distinct points, and the unique algebraic invariant of four points $u,v,x,y\in \Field\P^1$ is the cross-ratio:
\begin{equation}
\label{eq:bir} 
    \bir(u,v,x,y) 
    = \frac{(v-u)}{(v-x)}
    \div \frac{(y-u)}{(y-x)}
    \in \Field\P^1
\end{equation}
It satisfies in particular $\bir(z,0,1,\infty)=z$ and $\bir(z,0,w,\infty)=z/w$, whence the cocycle rule:
\begin{equation*}
    \bir(z,v,x,y)=\bir(z,v,w,y) \bir(w,v,x,y).
\end{equation*}

For a triple $(u,v,w)$ of points in $\Field\P^1$, we define their Maslov index in $\{0\}\sqcup\Field^\times/(\Field^\times)^2$ by lifting them to vectors $(\vec{u},\vec{v},\vec{w})$ in $\Field^2$ which are not all zero but with zero sum, and taking the determinant $\det(\vec{u},\vec{v})=\det(\vec{v},\vec{w})=\det(\vec{w},\vec{u})$.
It is preserved by the subgroup $\PSL_2(\Field)$, which acts freely transitively on the triples of distinct points with a given Maslov index.

A non-trivial $A\in \PSL_2(\Field)$ has two fixed points $\alpha_-,\alpha_+$ in the projective line over $\sqrt{\Field}$, that are well defined up to permutation.
The element $\epsilon_A^{\pm 2} := \bir(Ax,\alpha_\mp,x,\alpha_\pm)$ does not depend on $x\in \Field \P^1 \setminus\{\alpha_-,\alpha_+\}$, and is well defined up to inversion: it is called the \emph{period} of $A$.
The unique algebraic function on $\PSL_2(\Field)$ which is invariant by conjugacy is the discriminant:
\begin{equation*}
    \disc(A) = \left(\epsilon_A-\epsilon_A^{-1}\right)^2 
\end{equation*}
whose class in $\{0\}\sqcup \Field^\times /(\Field^\times)^2$ defines the \emph{type} of $A$.
All non-trivial elements with $\disc=0$ are conjugate. The elements with $\disc \ne 0$ are called \emph{semi-simple}.
The conjugacy classes of elements with $\disc \equiv 1 \bmod{(\Field^\times)^2}$ are uniquely characterised by the value of their discriminant.

Consider $A,B\in \PSL_2(\Field)$ of the same type, and fix a square root of $\disc(A)\disc(B)\in (\Field^\times)^2$.
Then one may order their fixed points $(\alpha_-,\alpha_+)$ and $(\beta_-,\beta_+)$ up to simultaneous inversion, and consistently define their cross-ratio $\bir(A,B)$ by:
\begin{equation*}
    \bir(A,B):= \bir(\alpha_+,\alpha_-,\beta_+,\beta_-) \in \Field\P^1
\end{equation*}
which is $\notin \{0,1,\infty\}$ unless $A$ and $B$ share a fixed point, and satisfies the symmetry property:
\begin{equation*}
    \frac{1}{\bir(A,B)}+\frac{1}{\bir(A,B^{-1})}=1.
\end{equation*}
We may also define their cosine (using their adjoint action on the Lie algebra $\Sl_2(\Field)$ as in \cite{CLS_Conj-PSL2K_2022}), which is related to the cross-ratio by:
\begin{equation}
\label{eq:cos-bir} 
    \cos(A,B)=\frac{1}{\bir(A,B)}-\frac{1}{\bir(A,B^{-1})}.
\end{equation}

\begin{Theorem}[\cite{CLS_Conj-PSL2K_2022}]
\label{Thm:conj-PSL_2(Z)}
Two semi-simple elements $A,B\in \PSL_2(\Field)$ are conjugate if and only if $\disc(A)=\Delta=\disc(A)$ and $\bir(A,B)\equiv 1 \bmod{\Norm_\Field \Field[\sqrt{\Delta}]}$.

More precisely, if $\Delta\ne 0$ and $\bir(A,B)\notin\{1,\infty\}$, then we may lift $A,B$ in $\SL_2(\Field)$ so that $\Aa=A-\tfrac{1}{2}\Tr(A), \Bb=B-\tfrac{1}{2}\Tr(A)$ satisfy $\Tr(\Aa\Bb)=-\Delta\cos(A,B)$, and the elements $C\in \SL_2(\Field)$ such that $CA C^{-1}= B$ are parametrized by the $\Field$-points of the Pell-Fermat conic:
\begin{equation*}
    (x,y)\in \Field\times \Field 
    \: \colon\: \quad
    (2x)^2+\Delta y^2 = \bir(A,B)
    \quad \text{according to} \quad
    C(x,y) = x(\Id+\Bb\Aa^{-1})+y(\Aa+\Bb).
\end{equation*}
\end{Theorem}

\subsection{Over the real field}

The automorphism group $\PGL_2(\C)\simeq \PSL_2(\C)$ of the complex projective line $\C\P^1$ contains $\PGL_2(\R)$ as the stabiliser of the real projective line $\R\P^1$.
The index-two subgroup $\PSL_2(\R)$ also preserves the upper half-plane $\H\P=\{z\in \C \mid \Im(z)>0\}\subset \C\P^1$, or equivalently the orientation induced on its boundary $\partial \H\P$, also given by the cyclic order $\cord(x,y,z) \in \{\pm 1\}$ of any triple of distinct points $x,y,z$ of $\R\P^1$ (which equals their Maslov index).

The complex structure on $\H\P$ is conformal to a unique hyperbolic metric. The hyperbolic distance $\lambda$ between $w,z\in \H\P$ can be deduced from the cross-ratio by $\bir(\Bar{z},z,\Bar{w},w)^{-1}=\left(\cosh \tfrac{\lambda }{2}\right)^{2}$.
This realises $\PSL_2(\R)$ as the positive isometry group of the hyperbolic plane: it preserves the previous cross-ratio and acts simply-transitively on positive triples of distinct points in $\R\P^1$, thus it preserves the hyperbolic metric and acts simply transitively on the unit tangent bundle of $\H\P$.

The type of $A \in \PSL_2(\R)$ is elliptic or parabolic or hyperbolic according to the value of $\sign \disc(A) \in \{-1,0,1\}$, equal to the number of distinct fixed points in $\R\P^1$ minus $1$.

A hyperbolic $A \in \PSL_2(\R)$ acts on $\H\P$ by translation along an (oriented) geodesic $\alpha$ whose endpoints $\alpha_-,\alpha_+ \in \R\P^1$ are its repulsive and attractive fixed points.
With this order the period satisfies $\epsilon_A^2 >1$.
The translation length $\lambda_A = \log(\epsilon_A^2)$ yields $\disc(A)=\left(2\sinh \tfrac{1}{2}\lambda_A\right)^2$.

\begin{Lemma}[Real geometry of the cross-ratio]
\label{Lem:cos-cosh-sinh}
Consider hyperbolic $A,B\in \PSL_2(\R)$ whose axes $\alpha, \beta$ have distinct endpoints. For any lifts $A,B\in \SL_2(\R)$ we have:
\begin{equation*}\textstyle
    \cos(A,B)= \sign(\Tr(A)\Tr(B)) \frac{\Tr(AB)-\Tr(AB^{-1})}{\sqrt{\disc(A)\disc(B)}}.
\end{equation*}

If $\alpha$ and $\beta$ intersect, then it is with an angle $\theta\in ]-\pi,\pi[$ oriented anti-clockwise from $\alpha$ to $\beta$ according to $\sign(\theta)=:\sign(A,B)$, and 
satisfying $\cos(\theta)=\cos(A,B)$, thus:
\begin{equation*}
    \tfrac{1}{\bir(A,B)}
    = \tfrac{1 + \cos(\theta)}{2}
    = \left(\cos \tfrac{\theta}{2}\right)^{2}
\end{equation*}
Note that $\theta\mapsto -\theta$ by transposition of $(A,B)$ whereas $\theta\mapsto \sign(\theta)\pi-\theta$ by inversion of $A$ or $B$.

If $\alpha$ and $\beta$ are disjoint, then they are joint by a unique ortho-geodesic arc with length $\lambda \in \R$,
whose $\sign(\lambda)=\sign \bir(A,B)$ compares the co-orientations induced by each axis, and satisfying $\cos(A,B)=\sign(\lambda) \cosh (\lambda)$, thus respectively:
\begin{equation*}
    \tfrac{1}{\bir(A,B)}
    = \tfrac{1 + \cosh(\lambda)}{2}
    = \left(\cosh \tfrac{\lambda}{2}\right)^{2}
    \quad \mathrm{and}\quad
    \tfrac{1}{\bir(A,B)}
    = \tfrac{1 - \cosh(\lambda)}{2}
    = -\left(\sinh \tfrac{\lambda}{2}\right)^{2}
\end{equation*}
Note that $\lambda$ is invariant by transposition of $(A,B)$ whereas $\lambda\mapsto -\lambda$ by inversion of $A$ or $B$.
\end{Lemma}
\begin{figure}[h]
    \centering
    \scalebox{.38}{
\begin{tikzpicture}[line cap=round,line join=round,>=triangle 45,x=4.0cm,y=4.0cm]
\clip(-1.2,-1.2) rectangle (1.3,1.3);
\draw [shift={(0.41,0.5)},line width=1.8pt,color=red,fill=red,fill opacity=0.1] (0,0) -- (-40.55:0.12) arc (-40.55:67.61:0.12) -- cycle;
\draw [shift={(0.41,0.5)},line width=1.8pt,color=red,fill=red,fill opacity=0.1] (0,0) -- (139.45:0.12) arc (139.45:247.61:0.12) -- cycle;
\draw [line width=1.2pt] (0,0) circle (4cm);
\draw [line width=2pt,color=blue, ->] (-0.15,0.99)-- (1,0);
\draw [line width=2pt,color=blue, ->] (-0.2,-0.98)-- (0.55,0.83);
\draw [line width=1.2pt,dash pattern=on 3pt off 3pt,color=blue] (-0.15,0.99)-- (0.55,0.83);
\draw [line width=1.2pt,dash pattern=on 3pt off 3pt,color=blue] (-0.2,-0.98)-- (1,0);
\begin{scriptsize}
\draw[color=blue] (-0.32,1.08) node {\Huge$\alpha_-$};
\draw[color=blue] (1.2,-0.03) node {\Huge$\alpha_+$};
\draw[color=blue] (-0.32,-1.1) node {\Huge$\beta_-$};
\draw[color=blue] (0.78,0.92) node {\Huge$\beta_+$};
\draw[color=red] (0.68,0.55) node {\Huge$\theta$};
\end{scriptsize}
\end{tikzpicture}
}
    \hfill
    \scalebox{.38}{
\begin{tikzpicture}[line cap=round,line join=round,>=triangle 45,x=4.0cm,y=4.0cm]
\clip(-1.2,-1.2) rectangle (1.3,1.3);
\draw [line width=1.2pt] (0,0) circle (4cm);
\draw [line width=1.2pt,dash pattern=on 3pt off 3pt,color=blue] (-0.15,0.99)-- (1,0);
\draw [line width=1.2pt,dash pattern=on 3pt off 3pt,color=blue] (-0.2,-0.98)-- (0.55,0.83);

\draw [line width=2pt,color=blue,->] (-0.15,0.99)-- (0.55,0.83);
\draw [line width=2pt,color=blue,->] (-0.2,-0.98)-- (1,0);

\draw [line width=2pt,color=red] (0.68,-0.27)-- (0.28,0.89);
\begin{scriptsize}
\draw[color=blue] (-0.32,1.08) node {\Huge$\alpha_-$};
\draw[color=blue] (1.19,0.05) node {\Huge$\beta_+$};
\draw[color=blue] (-0.32,-1.08) node {\Huge$\beta_-$};
\draw[color=blue] (0.87,0.84) node {\Huge$\alpha_+$};
\draw[color=red] (0.62,0.57) node {\Huge$\lambda$};
\end{scriptsize}
\end{tikzpicture}
}
    \scalebox{.38}{
\begin{tikzpicture}[line cap=round,line join=round,>=triangle 45,x=4.0cm,y=4.0cm]
\clip(-1.2,-1.2) rectangle (1.3,1.3);
\draw [line width=1.2pt] (0,0) circle (4cm);
\draw [line width=1.2pt,dash pattern=on 3pt off 3pt,color=blue] (-0.15,0.99)-- (1,0);
\draw [line width=1.2pt,dash pattern=on 3pt off 3pt,color=blue] (-0.2,-0.98)-- (0.55,0.83);

\draw [line width=2pt,color=blue,->] (-0.15,0.99)-- (0.55,0.83);
\draw [line width=2pt,color=blue,<-] (-0.2,-0.98)-- (1,0);

\draw [line width=2pt,color=red] (0.68,-0.27)-- (0.28,0.89);
\begin{scriptsize}
\draw[color=blue] (-0.32,1.08) node {\Huge$\alpha_-$};
\draw[color=blue] (1.19,0.05) node {\Huge$\beta_-$};
\draw[color=blue] (-0.32,-1.08) node {\Huge$\beta_+$};
\draw[color=blue] (0.87,0.84) node {\Huge$\alpha_+$};
\draw[color=red] (0.62,0.57) node {\Huge$\lambda$};
\end{scriptsize}
\end{tikzpicture}
}
    \caption*{Angle $\theta\in \,]0,\pi[$. Ortho-geodesics: co-oriented ($\lambda>0$) and disco-oriented ($\lambda<0$).}
\end{figure}
\begin{proof}
    The statements are invariant by the $\PGL_2(\R)$ action so we may fix three points among $(\alpha_-, \alpha_+,\beta_-, \beta_+)$ and compute in terms of the last.
    The cosine formula was used in \cite{Wolpert_formula-cosine-Fenchel-Nielsen_1982}.
\end{proof}

\section{The modular group \texorpdfstring{$\PSL_2(\Z)$}{PSL(2,Z)}}

\subsection{The modular orbifold}

The subgroup $\PSL_2(\Q)$ of $\PSL_2(\R)$ is the stabiliser of the rational projective line $\Q\P^1$.
The discrete subgroup $\PSL_2(\Z)$ is the stabiliser of the ideal triangulation $\Tri$ of $\H\P$ with vertex set $\Q\P^1$ and edges all geodesics whose endpoints $\tfrac{p}{q},\tfrac{r}{s}$ satisfy $\lvert ps-qr \rvert =1$.

Consider the action of $\PSL_2(\Z)$ on $\Tri$.
It is transitive on the set of edges, which is in bijection with the orbit of $i\in (0,\infty)$, and the stabiliser of $i$ is the subgroup of order $2$ generated by $S$.
It is transitive on the set of triangles, which is in bijection with the orbit of $j=\exp(i\pi/3)$, and the stabiliser of $j$ is the subgroup of order $3$ generated by $T$.
Thus it is freely transitive on the flags of $\Tri$, or equivalently on the oriented edges, and we deduce that $\PSL_2(\Z)=\Z/2*\Z/3$ is the free amalgam of its subgroups generated by $S$ and $T$.
\begin{equation*}
    S = \begin{pmatrix}
    0 & -1 \\ 1 & 0
    \end{pmatrix}
    \qquad 
    T = \begin{pmatrix}
    1 & -1 \\ 1 & 0
    \end{pmatrix}
\end{equation*}

We also find that $\PSL_2(\Z)$ acts properly discontinuously on $\H\P$ with fundamental domain the triangle $(\infty,0,j)$.
We may cut it along the geodesic arc $(i,j)$ to obtain a pair of isometric triangles $(i,j,\infty)$ and $(i,j,0)$. Identifying them along their isometric edges yields the \emph{modular orbifold}
\begin{equation*}
    \M = \PSL_2(\Z)\backslash \H\P.
\end{equation*}
It is a hyperbolic two-dimensional orbifold, with conical singularities of order $2$ and $3$ associated to the fixed points $i$ and $j$ of $S$ and $T$, and a cusp associated to the fixed point $\infty \in \partial \H\P$ of $R=TS^{-1}$.

\begin{figure}[h]
    \centering
    \scalebox{0.48}{
\begin{tikzpicture}[line cap=round,line join=round,>=triangle 45,x=4.166666666666667cm,y=4.166666666666667cm]
\clip(-1.3,-1.3) rectangle (1.3,1.3);
\draw [line width=1.5pt] (0.,0.) circle (4.166666666666667cm);
\draw [line width=2.pt,color=forestgreen] (0.,1.)-- (0.8,0.6);
\draw [line width=2.pt,color=forestgreen] (0.8,0.6)-- (1.,0.);
\draw [line width=2.pt,color=forestgreen] (1.,0.)-- (0.8,-0.6);
\draw [line width=2.pt,color=forestgreen] (0.8,-0.6)-- (0.,-1.);
\draw [line width=2.pt,color=forestgreen] (0.,-1.)-- (-0.8,-0.6);
\draw [line width=2.pt,color=forestgreen] (-1.,0.)-- (-0.8,-0.6);
\draw [line width=2.pt,color=forestgreen] (-0.8,0.6)-- (-1.,0.);
\draw [line width=2.pt,color=forestgreen] (-0.8,0.6)-- (0.,1.);
\draw [line width=2.pt,color=forestgreen] (0.,1.)-- (0.,-1.);
\draw [line width=2.pt,color=forestgreen] (-1.,0.)-- (0.,1.);
\draw [line width=2.pt,color=forestgreen] (-1.,0.)-- (0.,-1.);
\draw [line width=2.pt,color=forestgreen] (0.,-1.)-- (1.,0.);
\draw [line width=2.pt,color=forestgreen] (-0.8,0.6)-- (0.,1.);
\draw [line width=2.pt,color=forestgreen] (0.,1.)-- (0.8,0.6);
\draw [line width=2.pt,color=forestgreen] (0.8,0.6)-- (1.,0.);
\draw [line width=2.pt,color=forestgreen] (1.,0.)-- (0.8,-0.6);
\draw [line width=2.pt,color=forestgreen] (0.8,-0.6)-- (0.,-1.);
\draw [line width=2.pt,color=forestgreen] (0.,-1.)-- (-0.8,-0.6);
\draw [line width=2.pt,color=forestgreen] (-0.8,-0.6)-- (-1.,0.);
\draw [line width=2.pt,color=forestgreen] (-1.,0.)-- (-0.8,0.6);
\draw [line width=2.pt,color=forestgreen] (0.,1.)-- (1.,0.);
\draw [line width=2.pt,color=forestgreen] (-0.8,0.6)-- (-0.4705882352941182,0.8823529411764702);
\draw [line width=2.pt,color=forestgreen] (-0.4705882352941182,0.8823529411764702)-- (0.,1.);
\draw [line width=2.pt,color=forestgreen] (0.,1.)-- (0.47058823529411786,0.8823529411764705);
\draw [line width=2.pt,color=forestgreen] (0.47058823529411786,0.8823529411764705)-- (0.8,0.6);
\draw [line width=2.pt,color=forestgreen] (0.8,0.6)-- (0.9600317911655876,0.27989097868883483);
\draw [line width=2.pt,color=forestgreen] (0.9600317911655876,0.27989097868883483)-- (1.,0.);
\draw [line width=2.pt,color=forestgreen] (1.,0.)-- (0.96,-0.28);
\draw [line width=2.pt,color=forestgreen] (0.96,-0.28)-- (0.8,-0.6);
\draw [line width=2.pt,color=forestgreen] (0.8,-0.6)-- (0.47201199143714934,-0.8815921278797451);
\draw [line width=2.pt,color=forestgreen] (0.47201199143714934,-0.8815921278797451)-- (0.,-1.);
\draw [line width=2.pt,color=forestgreen] (0.,-1.)-- (-0.47201199143714934,-0.8815921278797451);
\draw [line width=2.pt,color=forestgreen] (-0.47201199143714934,-0.8815921278797451)-- (-0.8,-0.6);
\draw [line width=2.pt,color=forestgreen] (-0.96,-0.28)-- (-0.8,-0.6);
\draw [line width=2.pt,color=forestgreen] (-1.,0.)-- (-0.96,-0.28);
\draw [line width=2.pt,color=forestgreen] (-0.9600415291452317,0.27985757505646536)-- (-1.,0.);
\draw [line width=2.pt,color=forestgreen] (-0.8,0.6)-- (-0.9600415291452317,0.27985757505646536);
\begin{scriptsize}
\draw [fill=black] (0.8,-0.6) circle (2.5pt);

\draw [fill=black] (0.47058823529411786,0.8823529411764705) circle (2.5pt);
\draw[color=black] (0.52,1) node {\huge$\frac{3}{1}$};

\draw [fill=black] (0.47201199143714934,-0.8815921278797451) circle (2.5pt);
\draw[color=black] (0.52,-1) node {\huge$\frac{1}{3}$};

\draw [fill=black] (-0.47201199143714934,-0.8815921278797451) circle (2.5pt);
\draw[color=black] (-0.57,-1) node {\huge$-\frac{1}{3}$};

\draw [fill=black] (-0.4705882352941182,0.8823529411764702) circle (2.5pt);
\draw[color=black] (-0.57,1) node {\huge$-\frac{3}{1}$};

\draw [fill=black] (0.9600317911655876,0.27989097868883483) circle (2.5pt);
\draw[color=black] (1.03,0.35) node {\huge$\frac{3}{2}$};

\draw [fill=black] (0.96,-0.28) circle (2.5pt);
\draw[color=black] (1.03,-0.35) node {\huge$\frac{2}{3}$};

\draw [fill=black] (-0.9600415291452317,0.27985757505646536) circle (2.5pt);
\draw[color=black] (-1.1,0.35) node {\huge$-\frac{3}{2}$};

\draw [fill=black] (-0.96,-0.28) circle (2.5pt);
\draw[color=black] (-1.1,-0.35) node {\huge$-\frac{2}{3}$};

\draw [fill=black] (-0.8,0.6) circle (2.5pt);
\draw[color=black] (-0.95,0.7) node {\huge$-\frac{2}{1}$};
\draw [fill=black] (0.,1.) circle (2.5pt);
\draw[color=black] (0.,1.15) node {\huge$\frac{1}{0}$};
\draw [fill=black] (1.,0.) circle (2.5pt);
\draw[color=black] (1.1,0.) node {\huge$\frac{1}{1}$};
\draw [fill=black] (0.,-1.) circle (2.5pt);
\draw[color=black] (0.,-1.15) node {\huge$\frac{0}{1}$};
\draw [fill=black] (-1,0.) circle (2.5pt);
\draw[color=black] (-1.15,0.) node {\huge$-\frac{1}{1}$};
\draw [fill=black] (0.8,0.6) circle (2.5pt);
\draw[color=black] (0.9,0.65) node {\huge$\frac{2}{1}$};
\draw [fill=black] (-0.8,-0.6) circle (2.5pt);
\draw[color=black] (-0.95,-0.7) node {\huge$-\frac{1}{2}$};
\draw [fill=black] (0.8,-0.6) circle (2.5pt);
\draw[color=black] (0.9,-0.65) node {\huge$\frac{1}{2}$};
\end{scriptsize}
\end{tikzpicture}
}
    \hfill
    \scalebox{0.48}{\input{images/tikz/PSL2Z-pavage-PH-fundom_b}}
    \caption*{The ideal triangulation of $\H\P$ together with its dual trivalent tree $\Tree$ yield the modular tessellation with fundamental domain $(0,j,\infty)$.}
    \label{fig:LagranTree}
\end{figure}

The modular group $\PSL_2(\Z)$ is the orbifold fundamental group of $\M$, so its conjugacy classes correspond to the free homotopy classes of (oriented) loops in $\M$.
The elliptic conjugacy classes are those of $S$ and $T^{\pm 1}$ which correspond to loops encircling the singularities, and the parabolic conjugacy classes are those of $R^n$ which correspond to loops encircling the cusp.
The conjugacy class of a hyperbolic $A\in \PSL_2(\Z)$ corresponds to the free homotopy class of a unique closed geodesic $[\alpha]\subset \M$, and its length equals $\lambda_A=\log(\epsilon_A^2)$. These are the so called \emph{modular geodesics}.

\subsection{Acting on a trivalent tree}

The preimage of the segment $(i,j)\subset \M$ in $\H\P$ forms a bipartite tree $\Tree'$, the first barycentric subdivision of a trivalent tree $\Tree$ which is dual to the ideal triangulation.
The \emph{base edge} $(i,j)$ of $\Tree'$ defines the \emph{oriented base edge} $\vec{e}_i$ of $\Tree$.
The action of $\PSL_2(\Z)$ is freely transitive on the set of edges of $\Tree'$ hence on the set of oriented edges of $\Tree$.
It also preserves the cyclic order on the set of edges incident to each vertex (given by the surface embedding $\Tree \subset \H\P$).
This is equivalent to the cyclic order function $\cord(x,y,z)\in \{-1,0,1\}$ of three points $x,y,z\in \Tree \cup \partial \Tree$.
Thus $\PSL_2(\Z)$ is the full automorphism group of cyclically ordered simplicial tree $(\Tree, \cord)$. 

We may now use this action to find some conjugacy invariants of primitive elements by considering their stable subsets.
Recall that an element in $\PSL_2(\Z)$ is called primitive when it generates a maximal cyclic subgroup.
The (primitive) elliptic conjugacy classes correspond to the vertices of $\Tree'$ and the primitive parabolic conjugacy classes correspond to the connected components of $\H\P\setminus \Tree$.

Let $A\in \PSL_2(\Z)$ be primitive of infinite order.
It acts on $\Tree$ by translation along an oriented geodesic $\alpha_\Tree$ called its \emph{combinatorial axis}, with endpoints $\alpha_-,\alpha_+\in \partial \Tree = \R\P^1$.
Observe that $\alpha_\Tree$ passes through the oriented base edge $\vec{e}_i$ of $\Tree$ exactly when its endpoints satisfy $\alpha_-\le 0 \le \alpha_+$, which is equivalent to saying that $A$ maps the base triangle $(0,1,\infty)$ to a triangle of the form $(\tfrac{b}{d}, \tfrac{a+b}{c+d}, \tfrac{a}{c})$ with $a,b,c,d\in \N$, in other terms that $A$ belongs to the monoid $\PSL_2(\N)$ freely generated by $\{L,R\}$.
In that case, $\alpha_\Tree$ follows a periodic sequence of left and right turns given by the $\{L,R\}$-factorisation of $A$, or the continued fraction expansion of the periodic number $\alpha_+$. 

The conjugacy class of $A$ corresponds to the orbit of $\alpha_\Tree$ under the action of $\PSL_2(\Z)$ on $\Tree$.
Hence the conjugacy classes of non-elliptic elements in $\PSL_2(\Z)$ correspond to the cyclic words over the alphabet $\{L,R\}$, and the hyperbolic classes yield the cycles in which both letters appear.
The linear representatives of such an $\{L,R\}$-cycle parametrize the intersection of the corresponding conjugacy class with $\PSL_2(\N)$, whose elements are called its \emph{Euclidean representatives}.
For infinite order $A\in \PSL_2(\Z)$, the minimum displacement length $d(x,A\cdot x)$ of a vertex $x\in \Tree$ equals the combinatorial length $\len(A) = \#R+\#L$ of a Euclidean representative.

Since the combinatorial and geometric axes of a hyperbolic $A\in \PSL_2(\Z)$ have the same endpoints $\alpha_-,\alpha_+ \in \R\P^1 = \partial \H\P = \partial \Tree$, they intersect the ideal triangulation in the same pattern. However the geometric axis also contains the information of the intersection pattern with its first barycentric subdivision, which is equivalent to the isotopy class of the modular geodesic in $\M$. 

\begin{figure}[h]
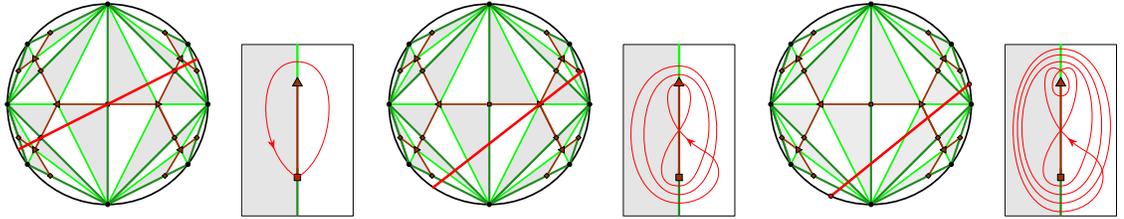

    \centering
    \scalebox{0.32}{\input{images/tikz/axis-RL-PH}}
    \scalebox{0.57}{
\begin{tikzpicture}[line cap=round,line join=round,>=triangle 45,rotate=90,x=2.0cm,y=2.0cm]
\clip(-1.1,-0.68) rectangle (1.1,0.68);

\draw [color=black] (1,0.65) -- (-1,.65) -- (-1,-0.65) -- (1,-0.65) -- cycle ;

\fill[line width=0.pt,color=marron,fill=grey,fill opacity=0.2] (1,0.65) -- (-1,.65) -- (-1,0.) -- (1,-0.) -- cycle ;

\draw[color=forestgreen,line width=1.5pt](-1,0) -- (-0. 45,0);
\draw[color=marron,line width=1.5pt](-0. 55,0) -- (0.  55,0);
\draw[color=green,line width=1.5pt](0.  55,0) -- (1,0);
\draw [fill=marron, rotate around={ 45:(-0.55,0)}] (-0.55,0) ++(-3.pt,0 pt) -- ++(3.pt,3.pt)--++(3.pt,-3.pt)--++(-3.pt,-3.pt)--++(-3.pt,3.pt);
\draw [fill=marron, rotate around={30:(0.  55,0))}] (0.  55,0) ++(0 pt,3.75pt) -- ++(3.2475952641916446pt,-5.625pt)--++(-6.495190528383289pt,0 pt) -- ++(3.2475952641916446pt,5.625pt);

\draw [color=red] (-0.55,0) to[out=-60, in=-90,looseness=1] (0.8,0)    to[out=90,in=60,looseness=1] (-0.55,0);

\draw[{Stealth[length=2mm,width=1.5mm]}-,color=red, shift={(-0.2,.278)},rotate=20] (-1pt,0) -- (1pt,0);


\end{tikzpicture}
}
    \hfill
    \scalebox{0.32}{\input{images/tikz/axis-RLL-PH}}
    \scalebox{0.57}{
\begin{tikzpicture}[line cap=round,line join=round,>=triangle 45,rotate=90,x=2.0cm,y=2.0cm]
\clip(-1.1,-0.68) rectangle (1.1,0.68);

\draw [color=black] (1,0.65) -- (-1,.65) -- (-1,-0.65) -- (1,-0.65) -- cycle ;

\fill[line width=0.pt,color=marron,fill=grey,fill opacity=0.2] (1,0.65) -- (-1,.65) -- (-1,0.) -- (1,-0.) -- cycle ;

\draw[color=forestgreen,line width=1.5pt](-1,0) -- (-0. 45,0);
\draw[color=marron,line width=1.5pt](-0. 55,0) -- (0.  55,0);
\draw[color=green,line width=1.5pt](0.  55,0) -- (1,0);
\draw [fill=marron, rotate around={ 45:(-0.  55,0)}] (-0.  55,0) ++(-3.pt,0 pt) -- ++(3.pt,3.pt)--++(3.pt,-3.pt)--++(-3.pt,-3.pt)--++(-3.pt,3.pt);
\draw [fill=marron, rotate around={30:(0.  55,0))}] (0.  55,0) ++(0 pt,3.75pt) -- ++(3.2475952641916446pt,-5.625pt)--++(-6.495190528383289pt,0 pt) -- ++(3.2475952641916446pt,5.625pt);

\draw [color=red] (-0.65,0 ) to[out=90, in=140,looseness=1.3] (0,0 );
\draw [color=red] (0,0 ) to[out=-40, in=-90,looseness=1.] (0.55,0 )    to[out=90,in=40,looseness=1.] (0,0 );
\draw [color=red] (-0.25,-0.3 ) to[out=160, in=-90,looseness=0.9] (-0.65,0 );
\draw [color=red] (-0.25,-0.3 ) to[out=-20, in=-90,looseness=0.9] (0.65,0 );
\draw [color=red] (0.65,0 ) to[out=90, in=90,looseness=1.] (-0.75,0 );
\draw [color=red] (-0.75,0 ) to[out=-90, in=-90,looseness=1] (0.75,0 );
\draw [color=red] (0.75,0 ) to[out=90, in=90,looseness=1.2] (-0.85,0 );
\draw [color=red] (-0.85,0 ) to[out=-90, in=-120,looseness=1.5] (-0.25,-0.3 );
\draw [color=red] (-0.25,-0.3) to[out=60, in=-140,looseness=1] (0,0 );

\draw[-{Stealth[length=2mm,width=1.5mm]}, color=red, shift={(-.107,-.1)},rotate=49] (-1pt,0) -- (1pt,0);


\end{tikzpicture}
}
    \hfill
    \scalebox{0.32}{\input{images/tikz/axis-RLLL-PH}}
    \scalebox{0.57}{
\begin{tikzpicture}[line cap=round,line join=round,>=triangle 45,rotate=90,x=2.0cm,y=2.0cm]
\clip(-1.1,-0.68) rectangle (1.1,0.68);

\draw [color=black] (1,0.65) -- (-1,.65) -- (-1,-0.65) -- (1,-0.65) -- cycle ;

\fill[line width=0.pt,color=marron,fill=grey,fill opacity=0.2] (1,0.65) -- (-1,.65) -- (-1,0.) -- (1,-0.) -- cycle ;

\draw[color=forestgreen,line width=1.5pt](-1,0) -- (-0. 45,0);
\draw[color=marron,line width=1.5pt](-0. 55,0) -- (0.  55,0);
\draw[color=green,line width=1.5pt](0.  55,0) -- (1,0);
\draw [fill=marron, rotate around={ 45:(-0.  55,0)}] (-0.  55,0) ++(-3.pt,0 pt) -- ++(3.pt,3.pt)--++(3.pt,-3.pt)--++(-3.pt,-3.pt)--++(-3.pt,3.pt);
\draw [fill=marron, rotate around={30:(0.  55,0))}] (0.  55,0) ++(0 pt,3.75pt) -- ++(3.2475952641916446pt,-5.625pt)--++(-6.495190528383289pt,0 pt) -- ++(3.2475952641916446pt,5.625pt);

\draw [color=red] (-0.65,0 ) to[out=90, in=140,looseness=1.3] (0,0 );
\draw [color=red] (0,0 ) to[out=-40, in=-60,looseness=1.2] (0.7,0 )  to[out=120, in=90,looseness=1.2] (0.4,0 ) to[out=-90, in=-120,looseness=1.2] (0.7,0 )  to[out=60,in=40,looseness=1.2] (0,0 );
\draw [color=red] (-0.25,-0.3 ) to[out=160, in=-90,looseness=0.9] (-0.65,0 );
\draw [color=red] (-0.25,-0.3 ) to[out=-20, in=-90,looseness=0.9] (0.8,0 );
\draw [color=red] (0.8,0 ) to[out=90, in=90,looseness=1.] (-0.75,0 );
\draw [color=red] (-0.75,0 ) to[out=-90, in=-90,looseness=1] (0.88,0 );
\draw [color=red] (0.88,0 ) to[out=90, in=90,looseness=1.] (-0.85,0 );
\draw [color=red] (-0.85,0 ) to[out=-90, in=-90,looseness=1.1] (0.95,0 );
\draw [color=red] (0.95,0 ) to[out=90, in=90,looseness=1.] (-0.95,0 );
\draw [color=red] (-0.95,0 ) to[out=-90, in=-120,looseness=1.5] (-0.25,-0.3 );
\draw [color=red] (-0.25,-0.3) to[out=60, in=-140,looseness=1] (0,0 );

\draw[-{Stealth[length=2mm,width=1.5mm]}, color=red, shift={(-.107,-.1)},rotate=49] (-1pt,0) -- (1pt,0);

\end{tikzpicture}
}
    \caption*{The geometric axes in $\P\H$ and their projections in $\M$ of $RL$, $RLL$, $RLLL$.}
\end{figure}

In \cite[Chapter 3]{CLS_phdthesis_2022}, we recover the isotopy class of $[\alpha]\subset \M$ from the $\{L,R\}$-cycle of $A$.
We also describe when $[\alpha]$ passes through the singular points $i$ or $j$.
Moreover \cite[Lemma 2.27]{CLS_phdthesis_2022} shows that if a (primitive) hyperbolic conjugacy class in $\PSL_2(\Z)$ is stable under inversion, then it contains (exactly $4$) symmetric matrices, and those all lie in $\PSL_2(\N)$ up to inversion.

\section{From the geometric cosine to the combinatorial cosign}

\subsection{The functions $\cross$ and $\cosign$.}

We now use the representation $\PSL_2(\Z)=\Aut(\Tree, \cord)$ to find conjugacy invariants for pairs of primitive infinite order elements by comparing the relative positions of their stable subsets, namely their combinatorial axes.

Let us first derive from the cyclic order function of three points, the crossing function of four points $u,v,x,y\in \Tree \cup \partial \Tree$ by:
\begin{equation}
\label{eq:cross} 
    \cross(u,v,x,y) =  \tfrac{1}{2}\left(\cord(u,x,v) - \cord(u,y,v)
\right)\end{equation} 
that is the algebraic intersection number of the oriented geodesics $(u,v)$ and $(x,y)$.
We denote $\across(u,v,x,y) \in \{0,\tfrac{1}{2},1\}$ the absolute value of $\cross(u,v,x,y)$ which yields the linking number of the cycles $(u,v)$, $(x,y)$ in the cyclically ordered boundary $\partial \Tree$.
One may compare the formula \eqref{eq:cross} defining $\cross$ with the formula \eqref{eq:bir} defining $\bir$ noticing that $\cord(x,y,z)=\sign \tfrac{y-z}{y-x}$. In particular, for $u,v,x,y\in \R\P^1$ we have \[\across(u,v,x,y) = 1 \iff \bir(u,v,x,y) > 1\]

Now consider two oriented bi-infinite geodesics $\alpha_\Tree$ and $\beta_\Tree$ of $\Tree$. Their intersection is either empty in which case we define $\cosign(\alpha_\Tree,\beta_\Tree)=0$, or else it consists in a geodesic containing at least one edge along which we may thus compare their orientations by $\cosign(\alpha_\Tree,\beta_\Tree)\in \{-1,+1\}$.

The functions $\cross$ and $\cosign$ are $\PSL_2(\Z)$-invariant, symmetric, and inverting the orientation of one argument results in a change of sign.

\begin{figure}[h]
    \centering
    \scalebox{.65}{
\begin{tabular}{
|>{\centering\arraybackslash}m{.2\textwidth}
|>{\centering\arraybackslash}m{.4\textwidth}
|>{\centering\arraybackslash}m{.4\textwidth}
|>{\centering\arraybackslash}m{.4\textwidth}
|}
\hline
\backslashbox[.2\textwidth]{$\cosign$}{$\cross$}
&$\mathbf{+1}$
&$\mathbf{0}$
&$\mathbf{-1}$\\
\hline
$\mathbf{+1}$ 
&\tikzset{every picture/.style={line width=0.75pt}} 

\begin{tikzpicture}[x=0.75pt,y=0.75pt,yscale=-1,xscale=1]

\draw [color=marron  ,draw opacity=1 ][line width=1.5]    (108.5,66.08) -- (198.5,66.08) ;
\draw [color=marron  ,draw opacity=1 ][line width=1.5]    (138.83,66.08) -- (138.83,51.25) ;
\draw [color=marron  ,draw opacity=1 ][line width=1.5]    (123.5,36.58) -- (138.83,51.25) ;
\draw [color=marron  ,draw opacity=1 ][line width=1.5]    (138.83,51.25) -- (153.5,36.58) ;
\draw [color=marron  ,draw opacity=1 ][line width=1.5]    (168.5,81.58) -- (168.33,66.08) ;
\draw [color=marron  ,draw opacity=1 ][line width=1.5]    (183.5,97.08) -- (168.5,81.58) ;
\draw [color=marron  ,draw opacity=1 ][line width=1.5]    (154,97.08) -- (168.5,81.58) ;
\draw [color=marron  ,draw opacity=1 ][line width=1.5]    (79,97.08) -- (109,66.08) ;
\draw [color=marron  ,draw opacity=1 ][line width=1.5]    (79,36.58) -- (109,66.08) ;
\draw [color=marron  ,draw opacity=1 ][line width=1.5]    (198.5,66.08) -- (228.5,36.58) ;
\draw [color=marron  ,draw opacity=1 ][line width=1.5]    (198.5,66.08) -- (228.5,97.08) ;

\draw [color=red  ,draw opacity=1 ][line width=2.25]    (83.5,33.58) .. controls (110.5,60.08) and (106,61.58) .. (129,60.58) .. controls (152,59.58) and (148,60.58) .. (153.5,66.08) .. controls (159,71.58) and (165.5,70.58) .. (180.5,70.58) .. controls (195.05,70.58) and (193.6,73.72) .. (219.97,99.01) ;
\draw [shift={(222.5,101.42)}, rotate = 223.48] [color=red  ,draw opacity=1 ][line width=2.25]    (17.49,-5.26) .. controls (11.12,-2.23) and (5.29,-0.48) .. (0,0) .. controls (5.29,0.48) and (11.12,2.23) .. (17.49,5.26)   ;
\draw [color=blue  ,draw opacity=1 ][line width=2.25]    (84.23,100.43) .. controls (110.8,72.73) and (106.31,71.23) .. (129.19,72.03) .. controls (152.08,72.82) and (148.09,71.83) .. (153.5,66.08) .. controls (158.91,60.34) and (165.39,61.3) .. (180.31,61.15) .. controls (194.7,61) and (196.98,57.15) .. (219.17,33.89) ;
\draw [shift={(221.67,31.28)}, rotate = 133.8] [color=blue  ,draw opacity=1 ][line width=2.25]    (17.49,-5.26) .. controls (11.12,-2.23) and (5.29,-0.48) .. (0,0) .. controls (5.29,0.48) and (11.12,2.23) .. (17.49,5.26)   ;
\draw  [color=violet  ,draw opacity=1 ][fill=violet  ,fill opacity=1 ] (152,65.33) .. controls (152.41,64.5) and (153.42,64.16) .. (154.25,64.58) .. controls (155.08,65) and (155.42,66) .. (155,66.84) .. controls (154.59,67.67) and (153.58,68) .. (152.75,67.59) .. controls (151.92,67.17) and (151.58,66.16) .. (152,65.33) -- cycle ;

\draw (69,16) node [anchor=north west][inner sep=0.75pt]  [color=red  ,opacity=1 ] [align=left] {$\displaystyle \alpha_-$};
\draw (69,102) node [anchor=north west][inner sep=0.75pt]  [color=blue  ,opacity=1 ] [align=left] {$\beta_+$};
\draw (224.5,16) node [anchor=north west][inner sep=0.75pt]  [color=blue  ,opacity=1 ] [align=left] {$\displaystyle \mathbf{\beta_+}$};
\draw (224,99.5) node [anchor=north west][inner sep=0.75pt]  [color=red  ,opacity=1 ] [align=left] {$\displaystyle \mathbf{\alpha_+}$};

\end{tikzpicture}
&

\tikzset{every picture/.style={line width=0.75pt}} 

\begin{tikzpicture}[x=0.75pt,y=0.75pt,yscale=-1,xscale=1]

\draw [color=marron  ,draw opacity=1 ][line width=1.5]    (108.5,66.08) -- (198.5,66.08) ;
\draw [color=marron  ,draw opacity=1 ][line width=1.5]    (138.83,66.08) -- (138.83,51.25) ;
\draw [color=marron  ,draw opacity=1 ][line width=1.5]    (123.5,36.58) -- (138.83,51.25) ;
\draw [color=marron  ,draw opacity=1 ][line width=1.5]    (138.83,51.25) -- (153.5,36.58) ;
\draw [color=marron  ,draw opacity=1 ][line width=1.5]    (168.5,81.58) -- (168.33,66.08) ;
\draw [color=marron  ,draw opacity=1 ][line width=1.5]    (183.5,97.08) -- (168.5,81.58) ;
\draw [color=marron  ,draw opacity=1 ][line width=1.5]    (154,97.08) -- (168.5,81.58) ;
\draw [color=marron  ,draw opacity=1 ][line width=1.5]    (79,97.08) -- (109,66.08) ;
\draw [color=marron  ,draw opacity=1 ][line width=1.5]    (79,36.58) -- (109,66.08) ;
\draw [color=marron  ,draw opacity=1 ][line width=1.5]    (198.5,66.08) -- (228.5,36.58) ;
\draw [color=marron  ,draw opacity=1 ][line width=1.5]    (198.5,66.08) -- (228.5,97.08) ;

\draw [color=blue  ,draw opacity=1 ][line width=2.25]    (84.5,32.08) .. controls (111.5,58.58) and (108.5,61.58) .. (131.5,60.58) .. controls (154.5,59.58) and (166.5,60.08) .. (181,60.58) .. controls (195.07,61.07) and (191.1,58.26) .. (220.05,34.36) ;
\draw [shift={(222.83,32.08)}, rotate = 140.91] [color=blue  ,draw opacity=1 ][line width=2.25]    (17.49,-5.26) .. controls (11.12,-2.23) and (5.29,-0.48) .. (0,0) .. controls (5.29,0.48) and (11.12,2.23) .. (17.49,5.26)   ;
\draw [color=red  ,draw opacity=1 ][line width=2.25]    (83.7,99.98) .. controls (110.79,73.48) and (107.78,70.48) .. (130.86,71.48) .. controls (153.94,72.48) and (165.99,71.98) .. (180.54,71.48) .. controls (194.65,71) and (190.97,74.37) .. (220.04,98.3) ;
\draw [shift={(222.83,100.58)}, rotate = 219] [color=red  ,draw opacity=1 ][line width=2.25]    (17.49,-5.26) .. controls (11.12,-2.23) and (5.29,-0.48) .. (0,0) .. controls (5.29,0.48) and (11.12,2.23) .. (17.49,5.26)   ;

\draw (69.6,12.1) node [anchor=north west][inner sep=0.75pt]  [color=blue  ,opacity=1 ] [align=left] {$\displaystyle \mathbf{\beta_-}$};
\draw (69.1,104.03) node [anchor=north west][inner sep=0.75pt]  [color=red  ,opacity=1 ] [align=left] {$\displaystyle \mathbf{\alpha_-}$};
\draw (231.2,12.1) node [anchor=north west][inner sep=0.75pt]  [color=blue  ,opacity=1 ] [align=left] {$\displaystyle \mathbf{\beta_+}$};
\draw (231.2,104.03) node [anchor=north west][inner sep=0.75pt]  [color=red  ,opacity=1 ] [align=left] {$\displaystyle \mathbf{\alpha_+}$};

\end{tikzpicture}
&\tikzset{every picture/.style={line width=0.75pt}} 

\begin{tikzpicture}[x=0.75pt,y=0.75pt,yscale=-1,xscale=1]

\draw [color=marron  ,draw opacity=1 ][line width=1.5]    (108.5,66.08) -- (198.5,66.08) ;
\draw [color=marron  ,draw opacity=1 ][line width=1.5]    (138.83,66.08) -- (138.83,51.25) ;
\draw [color=marron  ,draw opacity=1 ][line width=1.5]    (123.5,36.58) -- (138.83,51.25) ;
\draw [color=marron  ,draw opacity=1 ][line width=1.5]    (138.83,51.25) -- (153.5,36.58) ;
\draw [color=marron  ,draw opacity=1 ][line width=1.5]    (168.5,81.58) -- (168.33,66.08) ;
\draw [color=marron  ,draw opacity=1 ][line width=1.5]    (183.5,97.08) -- (168.5,81.58) ;
\draw [color=marron  ,draw opacity=1 ][line width=1.5]    (154,97.08) -- (168.5,81.58) ;
\draw [color=marron  ,draw opacity=1 ][line width=1.5]    (79,97.08) -- (109,66.08) ;
\draw [color=marron  ,draw opacity=1 ][line width=1.5]    (79,36.58) -- (109,66.08) ;
\draw [color=marron  ,draw opacity=1 ][line width=1.5]    (198.5,66.08) -- (228.5,36.58) ;
\draw [color=marron  ,draw opacity=1 ][line width=1.5]    (198.5,66.08) -- (228.5,97.08) ;

\draw [color=blue  ,draw opacity=1 ][line width=2.25]    (83.27,33.6) .. controls (110.27,60.1) and (106,61.58) .. (129,60.58) .. controls (152,59.58) and (148,60.58) .. (153.5,66.08) .. controls (159,71.58) and (165.5,70.58) .. (180.5,70.58) .. controls (195.05,70.58) and (194.33,73.89) .. (220.74,99.19) ;
\draw [shift={(223.27,101.6)}, rotate = 223.48] [color=blue  ,draw opacity=1 ][line width=2.25]    (17.49,-5.26) .. controls (11.12,-2.23) and (5.29,-0.48) .. (0,0) .. controls (5.29,0.48) and (11.12,2.23) .. (17.49,5.26)   ;
\draw [color=red  ,draw opacity=1 ][line width=2.25]    (83.85,100.43) .. controls (110.57,72.73) and (106.05,71.23) .. (129.06,72.03) .. controls (152.07,72.82) and (148.06,71.83) .. (153.5,66.08) .. controls (158.94,60.34) and (165.45,61.3) .. (180.45,61.15) .. controls (194.92,61) and (196.98,56.71) .. (219.57,34.54) ;
\draw [shift={(222.11,32.06)}, rotate = 135.77] [color=red  ,draw opacity=1 ][line width=2.25]    (17.49,-5.26) .. controls (11.12,-2.23) and (5.29,-0.48) .. (0,0) .. controls (5.29,0.48) and (11.12,2.23) .. (17.49,5.26)   ;
\draw  [color=violet  ,draw opacity=1 ][fill=violet  ,fill opacity=1 ] (152,65.33) .. controls (152.41,64.5) and (153.42,64.16) .. (154.25,64.58) .. controls (155.08,65) and (155.42,66) .. (155,66.84) .. controls (154.59,67.67) and (153.58,68) .. (152.75,67.59) .. controls (151.92,67.17) and (151.58,66.16) .. (152,65.33) -- cycle ;

\draw (69,16) node [anchor=north west][inner sep=0.75pt]  [color=blue ,opacity=1 ] [align=left] {$\displaystyle \mathbf{\beta_-}$};
\draw (69,102) node [anchor=north west][inner sep=0.75pt]  [color=red  ,opacity=1 ] [align=left] {$\alpha_-$};
\draw (224.5,16) node [anchor=north west][inner sep=0.75pt]  [color=red  ,opacity=1 ] [align=left] {$\displaystyle \mathbf{\alpha_+}$};
\draw (224,100.83) node [anchor=north west][inner sep=0.75pt]  [color=blue ,opacity=1 ] [align=left] {$\displaystyle \mathbf{\beta_+}$};

\end{tikzpicture}\\
\hline
$\mathbf{-1}$
&\tikzset{every picture/.style={line width=0.75pt}} 

\begin{tikzpicture}[x=0.75pt,y=0.75pt,yscale=-1,xscale=1]

\draw [color=marron  ,draw opacity=1 ][line width=1.5]    (108.5,66.08) -- (198.5,66.08) ;
\draw [color=marron  ,draw opacity=1 ][line width=1.5]    (138.83,66.08) -- (138.83,51.25) ;
\draw [color=marron  ,draw opacity=1 ][line width=1.5]    (123.5,36.58) -- (138.83,51.25) ;
\draw [color=marron  ,draw opacity=1 ][line width=1.5]    (138.83,51.25) -- (153.5,36.58) ;
\draw [color=marron  ,draw opacity=1 ][line width=1.5]    (168.5,81.58) -- (168.33,66.08) ;
\draw [color=marron  ,draw opacity=1 ][line width=1.5]    (183.5,97.08) -- (168.5,81.58) ;
\draw [color=marron  ,draw opacity=1 ][line width=1.5]    (154,97.08) -- (168.5,81.58) ;
\draw [color=marron  ,draw opacity=1 ][line width=1.5]    (79,97.08) -- (109,66.08) ;
\draw [color=marron  ,draw opacity=1 ][line width=1.5]    (79,36.58) -- (109,66.08) ;
\draw [color=marron  ,draw opacity=1 ][line width=1.5]    (198.5,66.08) -- (228.5,36.58) ;
\draw [color=marron  ,draw opacity=1 ][line width=1.5]    (198.5,66.08) -- (228.5,97.08) ;

\draw [color=blue  ,draw opacity=1 ][line width=2.25]    (224.44,100.39) .. controls (197.44,73.89) and (199.5,69.83) .. (176.5,70.83) .. controls (153.5,71.83) and (157.5,70.83) .. (152,65.33) .. controls (146.5,59.83) and (140,60.83) .. (125,60.83) .. controls (110.45,60.83) and (113.57,59.32) .. (87.3,34.13) ;
\draw [shift={(84.78,31.72)}, rotate = 43.48] [color=blue  ,draw opacity=1 ][line width=2.25]    (17.49,-5.26) .. controls (11.12,-2.23) and (5.29,-0.48) .. (0,0) .. controls (5.29,0.48) and (11.12,2.23) .. (17.49,5.26)   ;
\draw [color=red  ,draw opacity=1 ][line width=2.25]    (83.85,100.43) .. controls (110.57,72.73) and (106.05,71.23) .. (129.06,72.03) .. controls (152.07,72.82) and (148.06,71.83) .. (153.5,66.08) .. controls (158.94,60.34) and (165.45,61.3) .. (180.45,61.15) .. controls (194.92,61) and (196.98,56.71) .. (219.57,34.54) ;
\draw [shift={(222.11,32.06)}, rotate = 135.77] [color=red  ,draw opacity=1 ][line width=2.25]    (17.49,-5.26) .. controls (11.12,-2.23) and (5.29,-0.48) .. (0,0) .. controls (5.29,0.48) and (11.12,2.23) .. (17.49,5.26)   ;
\draw  [color=violet  ,draw opacity=1 ][fill=violet  ,fill opacity=1 ] (152,65.33) .. controls (152.41,64.5) and (153.42,64.16) .. (154.25,64.58) .. controls (155.08,65) and (155.42,66) .. (155,66.84) .. controls (154.59,67.67) and (153.58,68) .. (152.75,67.59) .. controls (151.92,67.17) and (151.58,66.16) .. (152,65.33) -- cycle ;

\draw (69,16) node [anchor=north west][inner sep=0.75pt]  [color=blue ,opacity=1 ] [align=left] {$\displaystyle \mathbf{\beta_+}$};
\draw (69,102) node [anchor=north west][inner sep=0.75pt]  [color=red  ,opacity=1 ] [align=left] {$\alpha_-$};
\draw (224.5,16) node [anchor=north west][inner sep=0.75pt]  [color=red  ,opacity=1 ] [align=left] {$\displaystyle \mathbf{\alpha_+}$};
\draw (224,100.83) node [anchor=north west][inner sep=0.75pt]  [color=blue,opacity=1 ] [align=left] {$\displaystyle \mathbf{\beta_-}$};

\end{tikzpicture}
&

\tikzset{every picture/.style={line width=0.75pt}} 

\begin{tikzpicture}[x=0.75pt,y=0.75pt,yscale=-1,xscale=1]

\draw [color=marron  ,draw opacity=1 ][line width=1.5]    (108.5,66.08) -- (198.5,66.08) ;
\draw [color=marron  ,draw opacity=1 ][line width=1.5]    (138.83,66.08) -- (138.83,51.25) ;
\draw [color=marron  ,draw opacity=1 ][line width=1.5]    (123.5,36.58) -- (138.83,51.25) ;
\draw [color=marron  ,draw opacity=1 ][line width=1.5]    (138.83,51.25) -- (153.5,36.58) ;
\draw [color=marron  ,draw opacity=1 ][line width=1.5]    (168.5,81.58) -- (168.33,66.08) ;
\draw [color=marron  ,draw opacity=1 ][line width=1.5]    (183.5,97.08) -- (168.5,81.58) ;
\draw [color=marron  ,draw opacity=1 ][line width=1.5]    (154,97.08) -- (168.5,81.58) ;
\draw [color=marron  ,draw opacity=1 ][line width=1.5]    (79,97.08) -- (109,66.08) ;
\draw [color=marron  ,draw opacity=1 ][line width=1.5]    (79,36.58) -- (109,66.08) ;
\draw [color=marron  ,draw opacity=1 ][line width=1.5]    (198.5,66.08) -- (228.5,36.58) ;
\draw [color=marron  ,draw opacity=1 ][line width=1.5]    (198.5,66.08) -- (228.5,97.08) ;

\draw [color=red  ,draw opacity=1 ][line width=2.25]    (84.5,32.08) .. controls (111.5,58.58) and (108.5,61.58) .. (131.5,60.58) .. controls (154.5,59.58) and (166.5,60.08) .. (181,60.58) .. controls (195.07,61.07) and (191,58.81) .. (219.95,34.94) ;
\draw [shift={(222.73,32.67)}, rotate = 140.91] [color=red  ,draw opacity=1 ][line width=2.25]    (17.49,-5.26) .. controls (11.12,-2.23) and (5.29,-0.48) .. (0,0) .. controls (5.29,0.48) and (11.12,2.23) .. (17.49,5.26)   ;
\draw [color=blue  ,draw opacity=1 ][line width=2.25]    (223.53,101.47) .. controls (196.69,74.97) and (199.14,70.58) .. (176.27,71.58) .. controls (153.41,72.58) and (141.48,72.08) .. (127.06,71.58) .. controls (113.08,71.1) and (115.2,75.58) .. (86.3,99.58) ;
\draw [shift={(83.53,101.87)}, rotate = 320.74] [color=blue  ,draw opacity=1 ][line width=2.25]    (17.49,-5.26) .. controls (11.12,-2.23) and (5.29,-0.48) .. (0,0) .. controls (5.29,0.48) and (11.12,2.23) .. (17.49,5.26)   ;

\draw (69.6,12.1) node [anchor=north west][inner sep=0.75pt]  [color=red  ,opacity=1 ] [align=left] {$\displaystyle \mathbf{\alpha_-}$};
\draw (69.1,104.03) node [anchor=north west][inner sep=0.75pt]  [color=blue  ,opacity=1 ] [align=left] {$\displaystyle \mathbf{\beta_+}$};
\draw (231.2,17.1) node [anchor=north west][inner sep=0.75pt]  [color=red  ,opacity=1 ] [align=left] {$\displaystyle \mathbf{\alpha_+}$};
\draw (231.2,104.03) node [anchor=north west][inner sep=0.75pt]  [color=blue ,opacity=1 ] [align=left] {$\displaystyle \mathbf{\beta_-}$};

\end{tikzpicture}
&

\tikzset{every picture/.style={line width=0.75pt}} 

\begin{tikzpicture}[x=0.75pt,y=0.75pt,yscale=-1,xscale=1]

\draw [color=marron  ,draw opacity=1 ][line width=1.5]    (108.5,66.08) -- (198.5,66.08) ;
\draw [color=marron  ,draw opacity=1 ][line width=1.5]    (138.83,66.08) -- (138.83,51.25) ;
\draw [color=marron  ,draw opacity=1 ][line width=1.5]    (123.5,36.58) -- (138.83,51.25) ;
\draw [color=marron  ,draw opacity=1 ][line width=1.5]    (138.83,51.25) -- (153.5,36.58) ;
\draw [color=marron  ,draw opacity=1 ][line width=1.5]    (168.5,81.58) -- (168.33,66.08) ;
\draw [color=marron  ,draw opacity=1 ][line width=1.5]    (183.5,97.08) -- (168.5,81.58) ;
\draw [color=marron  ,draw opacity=1 ][line width=1.5]    (154,97.08) -- (168.5,81.58) ;
\draw [color=marron  ,draw opacity=1 ][line width=1.5]    (79,97.08) -- (109,66.08) ;
\draw [color=marron  ,draw opacity=1 ][line width=1.5]    (79,36.58) -- (109,66.08) ;
\draw [color=marron  ,draw opacity=1 ][line width=1.5]    (198.5,66.08) -- (228.5,36.58) ;
\draw [color=marron  ,draw opacity=1 ][line width=1.5]    (198.5,66.08) -- (228.5,97.08) ;

\draw [color=red  ,draw opacity=1 ][line width=2.25]    (83.5,33.58) .. controls (110.5,60.08) and (106,61.58) .. (129,60.58) .. controls (152,59.58) and (148,60.58) .. (153.5,66.08) .. controls (159,71.58) and (165.5,70.58) .. (180.5,70.58) .. controls (195.05,70.58) and (193.6,73.72) .. (219.97,99.01) ;
\draw [shift={(222.5,101.42)}, rotate = 223.48] [color=red  ,draw opacity=1 ][line width=2.25]    (17.49,-5.26) .. controls (11.12,-2.23) and (5.29,-0.48) .. (0,0) .. controls (5.29,0.48) and (11.12,2.23) .. (17.49,5.26)   ;
\draw [color=blue  ,draw opacity=1 ][line width=2.25]    (223.15,31.74) .. controls (196.58,59.43) and (201.07,60.94) .. (178.18,60.14) .. controls (155.3,59.34) and (159.29,60.34) .. (153.88,66.08) .. controls (148.46,71.83) and (141.99,70.86) .. (127.07,71.02) .. controls (112.67,71.17) and (109.61,76.42) .. (87.36,99.79) ;
\draw [shift={(84.87,102.4)}, rotate = 313.8] [color=blue  ,draw opacity=1 ][line width=2.25]    (17.49,-5.26) .. controls (11.12,-2.23) and (5.29,-0.48) .. (0,0) .. controls (5.29,0.48) and (11.12,2.23) .. (17.49,5.26)   ;
\draw  [color=violet ,draw opacity=1 ][fill=violet ,fill opacity=1 ] (152,65.33) .. controls (152.41,64.5) and (153.42,64.16) .. (154.25,64.58) .. controls (155.08,65) and (155.42,66) .. (155,66.84) .. controls (154.59,67.67) and (153.58,68) .. (152.75,67.59) .. controls (151.92,67.17) and (151.58,66.16) .. (152,65.33) -- cycle ;

\draw (69,16) node [anchor=north west][inner sep=0.75pt]  [color=red ,opacity=1 ] [align=left] {$\displaystyle \mathbf{\alpha_-}$};
\draw (69.5,99.5) node [anchor=north west][inner sep=0.75pt]  [color=blue  ,opacity=1 ] [align=left] {$\displaystyle \mathbf{\beta_+}$};
\draw (224.5,16) node [anchor=north west][inner sep=0.75pt]  [color=blue ,opacity=1 ] [align=left] {$\displaystyle \mathbf{\beta_-}$};
\draw (224,99.5) node [anchor=north west][inner sep=0.75pt]  [color=red  ,opacity=1 ] [align=left] {$\displaystyle \mathbf{\alpha_+}$};

\end{tikzpicture}\\
\hline
\end{tabular}

}
    \caption*{Configurations of axes: $\cross$ and $\cosign$. Note that $\cross\ne 0 \implies \cosign = \pm 1$.}
\end{figure}

For hyperbolic $A,B \in \PSL_2(\Z)$ with axes $\alpha_\Tree=(\alpha_-,\alpha_+)$ and $\beta_\Tree=(\beta_-,\beta_+)$ in $\Tree$, we denote by $\cross(A,B)=\cross(\alpha_+,\alpha_-,\beta_+,\beta_-)$ and $\cosign(A,B)=\cosign(\alpha_\Tree,\beta_\Tree)$.
Notice $\cosign(A,B)=1$ if and only if there exists $C\in \PSL_2(\Z)$ such that $CAC^{-1}, CBC^{-1} \in \PSL_2(\N)$, in which case the set of such $C$ corresponds to the edges in $\alpha_\Tree \cap \beta_\Tree \subset \Tree$. 

\begin{Proposition}
\label{Prop:cosign(A,B)=sign(len(AB)-len(A/B))}
For hyperbolic $A,B \in \PSL_2(\Z)$ such that $\alpha_\Tree \cap \beta_\Tree \ne \emptyset$, we have:
\begin{equation*}
    \cosign(A,B)= \sign\left(\len AB -\len AB^{-1} \right).
\end{equation*}
\end{Proposition}

\begin{proof}
This follows from \cite[Proposition 1.6]{Paulin_Gromov-R-trees_1989}, which was corrected by \cite{Conder-Paulin_Erratum-Gromov-R-trees_2020}, and one may also consult \cite[Proposition 2.44]{CLS_phdthesis_2022}.
Compare with the cosine formula in Lemma \ref{Lem:cos-cosh-sinh}.
\end{proof}

\subsection{Deforming the \texorpdfstring{$\PSL_2(\Z)$}{PSL(2;Z)}-action on \texorpdfstring{$\H\P$}{HP} to the \texorpdfstring{$\PSL_2(\Z)$}{PSL(2;Z)}-action on \texorpdfstring{$\Tree$}{T}}

Let us define a family of representations $\rho_q \colon \SL_2(\Z) \to \SL_2(\R)$ depending algebraically on the parameter $q\in \R^*$ and with integral coefficients.
The Euclidean algorithm implies that $\SL_2(\Z)$ is generated by $\{S,R\}$, whence by $\{S,T\}$, or $\{L,R\}$.
Fix $S_q=S$ and let $T_q$ be the conjugate of $T$ by $\exp \tfrac{1}{2}\log(q)
\begin{psmallmatrix}
1&0\\0&-1
\end{psmallmatrix}$. 
Given $A\in \SL_2(\Z)$, we deduce $A_q=\rho_q(A)$ from any $\{S,T\}$-factorisation by replacing $T\mapsto T_q$, eg:
\begin{equation*}
    R_q = 
    \begin{pmatrix}
    q & 1 \\
    0 & q^{-1}
    \end{pmatrix}
    \quad \mathrm{and} \quad
    L_q =
    \begin{pmatrix}
    q & 0 \\
    1 & q^{-1}
    \end{pmatrix}.
\end{equation*}

This descends to a representation $\Bar{\rho}_q \colon \PSL_2(\Z)\to \PSL_2(\R)$ which is faithful and discrete (because $\disc(R_q) = (q-q^{-1})^2 \ge 0$), 
and positive in the sense that $T_q$ is a $2\pi/3$-rotation of $\H\P$ in the positive direction.
Conversely, every such representation is conjugate to $\Bar{\rho}_q$ for a unique $q>0$.

We have therefore parametrized the Teichm\"uller space of $\PSL_2(\Z)$ by the real algebraic set $\R_+^*$.
This Teichm\"uller space corresponds to the set of hyperbolic metrics $\M_q = \Bar{\rho}_q(\PSL_2(\Z))\backslash\H\P$ on the oriented modular orbifold as a topological space.
Observe intuitively that when $q\to \infty$, the hyperbolic orbifold $\M_q$ has a convex core which retracts onto the long geodesic arc $(i,j_q)$ connecting the conical singularities.
It lifts in $\H\P$ to an $\epsilon$-neighbourhood of a trivalent tree $\Tree_q$ with $\epsilon= \Theta\left(1/q^2\right)$.
Since the hyperbolic geodesics of $\M_q$ remain in its convex core, their angles must tend to $0\bmod{\pi}$.

\begin{figure}[h]
    \centering
    \scalebox{0.9}{

\begin{tikzpicture}[x=4cm,y=4cm,z=2cm]

\path (0.,0.,0.) coordinate (O)
(-0.5,0.,0.) coordinate (A)
(-0.8,-0.3,0.) coordinate (A')
(0.5,0.,0.) coordinate (B)
(0.8,-0.3,0.) coordinate (B')
(-0.8,0.3,0.) coordinate (C)
(0.8,0.3,0.) coordinate (D);

\fill[color=grey, opacity=0.2] (B') to[out=150, in=-90,looseness=0.7]  (B) arc(0:-180:0.5 and 0.1) (A) to[out=-90, in=30,looseness=0.7]  (A') -- (B') --  cycle;
\draw [color=forestgreen,line width=2pt,line cap =round] (A') to[out=30, in=-90,looseness=0.7]  (A) to[out=90,in=-30,looseness=0.7]  (C);
\draw [color=green,line width=2pt,line cap =round] (B') to[out=150, in=-90,looseness=0.7]  (B) to[out=90,in=-150,looseness=0.7]  (D);

\draw[name path=ElliBBack,line width=2.pt,color=violet,dashed] (0.5,0,0) arc(0:180:0.5 and 0.1);
\draw[line width=2.pt,name path=ElliBFront,color=violet] (0.5,0,0) arc(0:-180:0.5 and 0.1);

\draw[name path=arbre,line width=2pt, color=marron] (A') -- (B') ;

\node[draw,inner sep=3pt, fill=marron] at (A') {};
\node[isosceles triangle,isosceles triangle apex angle=60,draw,inner sep=2.5pt, fill=marron] at (B') {};

\end{tikzpicture}
}
    \qquad \qquad
    \scalebox{0.45}{

\begin{tikzpicture}[x=4cm,y=4cm,line cap=round,line join=round]

\clip (0.,0.) circle (1);

\path (0.,0.) coordinate (O)
(-0.55,0.) coordinate (A1)
(0.55,0.) coordinate (B1)
(0.72,-0.3) coordinate (B2)
(0.72,0.3) coordinate (B3)
(0.94,0.45) coordinate (D1)
(-0.94,-0.45) coordinate (D2)
(-0.94,0.45) coordinate (D1')
(0.94,-0.45) coordinate (D2');

\coordinate (A3) at ($(A1)!0.7!(D1')$);
\coordinate (A2) at ($(A1)!0.7!(D2)$);
\coordinate (B3) at ($(B1)!0.7!(D1)$);
\coordinate (B2) at ($(B1)!0.7!(D2')$);

\coordinate (C1) at ($(A1)!0.5!(D1')$);
\coordinate (C2) at ($(A1)!0.5!(D2)$);
\coordinate (C3) at ($(B1)!0.5!(D1)$);
\coordinate (C4) at ($(B1)!0.5!(D2')$);

\fill[color=grey, opacity=0.3] (D2) rectangle (D1);

\draw[line width=2pt,color=violet] (D2) rectangle (D1);

\draw[line width=2.5pt,color=black] (0.,0.) circle (1);



\draw[line width=1.5pt, color=marron] (A1) -- (B1) ;
\draw[line width=1.5pt, color=marron] (A1) -- (A2) ;
\draw[line width=1.5pt, color=marron] (A1) -- (A3) ;
\draw[line width=1.5pt, color=marron] (B1) -- (B2) ;
\draw[line width=1.5pt, color=marron] (B1) -- (B3) ;

\node[draw,inner sep=4pt, fill=marron] at (O) {};
\node[isosceles triangle,isosceles triangle apex angle=60,draw,inner sep=1.5pt,rotate=-5, fill=marron] at (A2) {};
\node[isosceles triangle,isosceles triangle apex angle=60,draw,inner sep=1.5pt,rotate=5, fill=marron] at (A3) {};
\node[isosceles triangle,isosceles triangle apex angle=60,draw,inner sep=1.5pt,rotate=-175, fill=marron] at (B2) {};
\node[isosceles triangle,isosceles triangle apex angle=60,draw,inner sep=1.5pt,rotate=175, fill=marron] at (B3) {};
\node[isosceles triangle,isosceles triangle apex angle=60,draw,inner sep=2.5pt,rotate=180, fill=marron] at (A1) {};
\node[isosceles triangle,isosceles triangle apex angle=60,draw,inner sep=2.5pt, fill=marron] at (B1) {};

\node[draw,inner sep=2pt, fill=marron,rotate=45] at (C1) {};
\node[draw,inner sep=2pt, fill=marron,rotate=45] at (C2) {};
\node[draw,inner sep=2pt, fill=marron,rotate=45] at (C3) {};
\node[draw,inner sep=2pt, fill=marron,rotate=45] at (C4) {};

\end{tikzpicture}
}
    \caption*{The convex core of $\M_q$ lifts in $\H\P$ to an $\epsilon$-neighbourhood of $\Tree_q$ with $\epsilon= \Theta\left(1/q^2\right)$.}
\end{figure}

To make this intuition precise, the geometric invariants $\disc$ and $\cos$ of $A_q,B_q$ define algebraic functions of $q$ whose degrees recover the combinatorial invariants $2\len$ and $\cosign$ of $A,B$.
This should not surprise someone acquainted with compactifications of Teichm\"uller space by actions on trees or by valuations \cite{Otal_compactification-varietes-representations_2015, MS_Aut(CV)_2020}.
Here the unique boundary point $q=\infty$ corresponds to the action on $\Tree$ or to the valuation $-\deg_q$.

\begin{Lemma}[From \texorpdfstring{$\cos$}{cos} to \texorpdfstring{$\cosign$}{cosign}]
\label{Lem:cos_q-limit}
Consider hyperbolic $A,B\in \PSL_2(\Z)$ such that $\across(A,B)=1$.
For all $q>0$ the elements $A_q,B_q\in \PSL_2(\R)$ are hyperbolic, and their geometric axes intersect at an angle whose cosine is an algebraic function of $q$ with limit $\cos(A_q,B_q) \xrightarrow[q\to \infty]{} \cosign(A,B)$.
\end{Lemma}

\begin{proof}
Lemma \ref{Lem:cos-cosh-sinh} expresses the cosine of the angles between the geometric axes of $A$ and $B$ as:
\begin{equation*}\textstyle
    \cos(A_q,B_q)
    = \sign(\Tr(A_q)\Tr(B_q))
    \frac{\Tr(A_qB_q)-\Tr(A_qB_q^{-1})}{\sqrt{\disc(A_q)\disc(B_q)}}.
\end{equation*}
For all $C\in \SL_2(\Z)$ the Laurent polynomial $\Tr(C_q)$ is reciprocal of degree $\len(C)$.
To find the limit as $q\to \infty$ we must compute the degrees and dominant terms of the polynomials involved in this expression.
Now recall from Proposition \ref{Prop:cosign(A,B)=sign(len(AB)-len(A/B))} that for hyperbolic $A,B \in \PSL_2(\Z)$ whose fixed points are linked we have $\cosign(A,B)= \sign\left(\len AB-\len AB^{-1} \right)$.
This completes the proof.
\end{proof}

\section{The unit tangent bundle to the modular orbifold}

\subsection{Modular knots and links}

The Lie group $\PSL_2(\R)$ identifies with the unit tangent bundle of the hyperbolic plane $\H\P$.
Its lattice $\PSL_2(\Z)$ acts on the left with quotient $\U=\PSL_2(\Z)\backslash\PSL_2(\R)$ the unit tangent bundle of the modular orbifold $\M= \PSL_2(\Z) \backslash \H\P$.
The fundamental group of $\U$ is the preimage of $\PSL_2(\Z)$ in the universal cover of $\PSL_2(\R)$, given by the central extension:
\begin{equation*}
    \Id \to \Z \to \widetilde{\PSL}_2(\Z) \to \PSL_2(\Z) \to \Id
\end{equation*}
and we find that $\pi_1(\U)$ is isomorphic to the braid group on three strands, hence to the fundamental group of a trefoil knot's complement.
In fact, the structure of the Seifert fibration $\U \to \M$ reveals that $\U$ is homeomorphic to the complement of a trefoil knot in the sphere (see \cite{Montesinos_Tesselations_1987, Dehornoy-Pinsky_template-pqr_2018} for such a proof).
In particular, any two disjoint loops in $\U$ have a well defined linking number.

The closed hyperbolic geodesics in $\M$ lift to the periodic orbits for the geodesic flow in its unit tangent bundle $\U$, and the primitive ones trace the so called \emph{modular knots}.
Together, they form the \emph{master modular link} whose components are indexed by the primitive hyperbolic conjugacy classes of the modular group.
We wish to relate the geometry and topology of the master modular link with the arithmetic and combinatorial properties of the modular group.

\begin{figure}[h]
    \centering
    \includegraphics[width=0.36\textwidth]{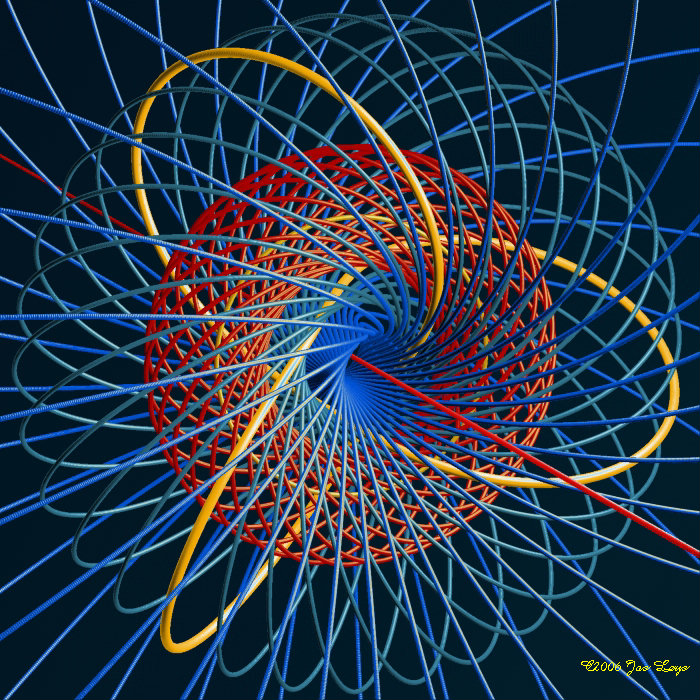}
    \includegraphics[width=0.48\textwidth]{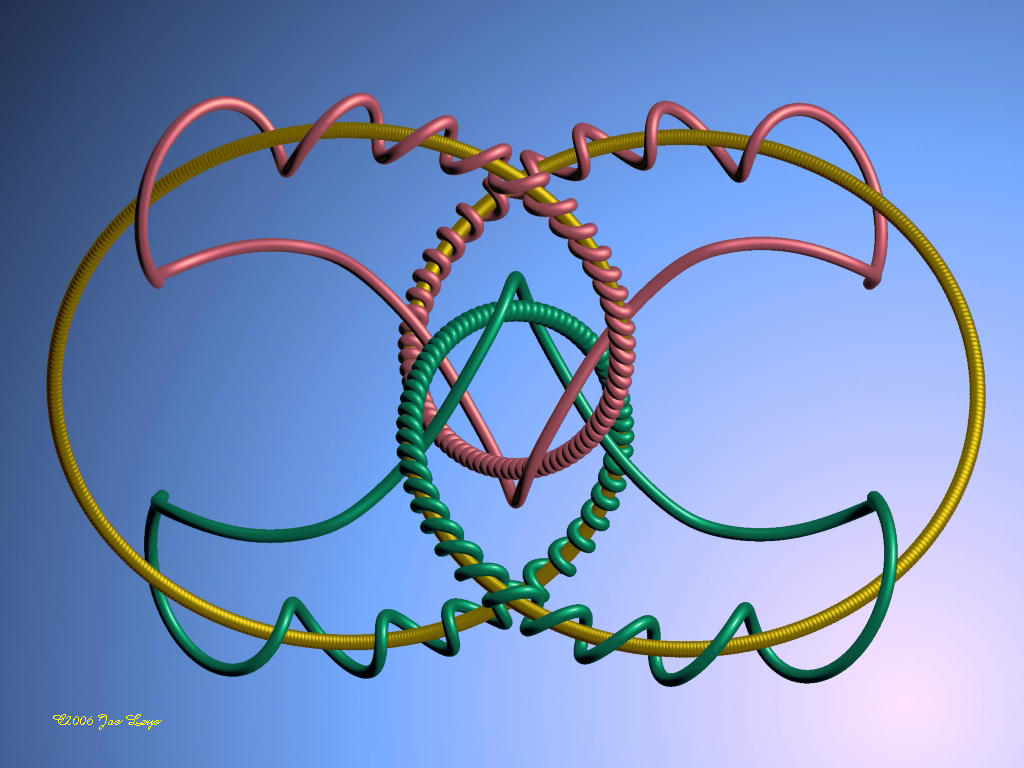}
    \caption*{The Seifert fibration $\U\to \M$ and two modular knots, from the \href{http://www.josleys.com/articles/ams_article/Lorenz3.htm}{online article} \cite{GhyLey_Lorenz-Modular-visual_2016} which proposes an animated introduction to the topology and dynamics of $\U$.}
\end{figure}

\subsection{The Lorenz template}

To describe the isotopy class of the master modular link, we rely on the construction of the Lorenz template and its embedding in $\U$, following \cite[§3.4]{Ghys_knots-dynamics_2006}.

The Lorenz template $\Lorenz$ is the branched surface obtained from the ideal triangle $(0,1,\infty)$ of $\H\P$ by identifying the side $(1,\infty)$ with the side $(0,\infty)$ through $R^{-1}$ and the side $(0,1)$ with the side $(0,\infty)$ through $L^{-1}$.
It is endowed with a semi-flow defined by the horizontal vector field whose periodic orbits correspond to the non-empty cycles on $\{L,R\}$ (this is an interval exchange map).
After its embedding $\Lorenz \hookrightarrow \U$ suggested in the following figure, those form \emph{the master Lorenz link}.

Consider a primitive hyperbolic conjugacy class in $\PSL_2(\Z)$: the geometric axes of its Euclidean representatives intersect the ideal triangle $(0,1,\infty)$ in a collection of segments which quotient to a closed connected loop in $\Lorenz$.
This loop is isotopic to the periodic orbit of the semi-flow indexed by the corresponding $\{L,R\}$-cycle.
More precisely, \'E. Ghys  \cite[§3.4]{Ghys_knots-dynamics_2006} showed the following.

\begin{Theorem}
The master modular link formed by all modular knots is isotopic to the master Lorenz link formed by the primitive periodic orbits of the semi-flow on the Lorenz template.

In particular, the Rademacher invariant of a primitive hyperbolic conjugacy class in $\PSL_2(\Z)$ equals the linking number between the corresponding modular knot and the trefoil.
\end{Theorem}

\begin{proof}[Outline of the proof]
The Fuchsian representation $\Bar{\rho}_q\colon \PSL_2(\Z) \to \PSL_2(\R)$ with quotient the hyperbolic orbifold $\M_q$ lifts to $\widetilde{\PSL}_2(\Z)\to \widetilde{\PSL}_2(\R)$ with quotient its unit tangent bundle $\U_q$. 
Varying $q\in ]1,+\infty[$ yields isotopies between the manifolds $\U_q$ which are all homeomorphic to the complement of a trefoil's neighbourhood, and conjugacies between their geodesic flows whose periodic orbits are indexed by the primitive conjugacy classes of infinite order in $\PSL_2(\Z)$.
As $q\to \infty$ the manifold $\U_q$ retracts onto a branched surface homeomorphic to the Lorenz template, and the master $q$-modular link isotopes to the periodic orbits of its semi-flow.
\end{proof}

\begin{figure}[t]
    \centering
    \includegraphics[width=0.48\textwidth]{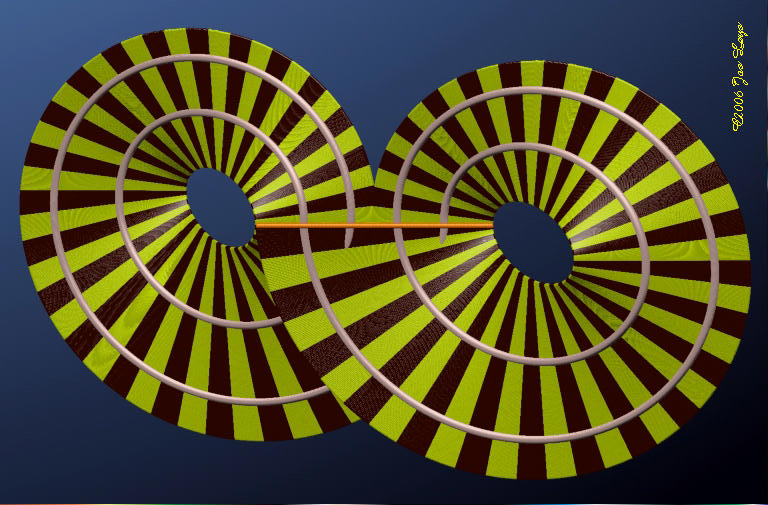}
    \caption*{Standard projection on $\S^2$ of the Lorenz template embedded in $\S^3$.}
    \label{fig:Lorenz-Template}
\end{figure}

\begin{figure}[h]
    \centering
    \includegraphics[width=0.42\textwidth]{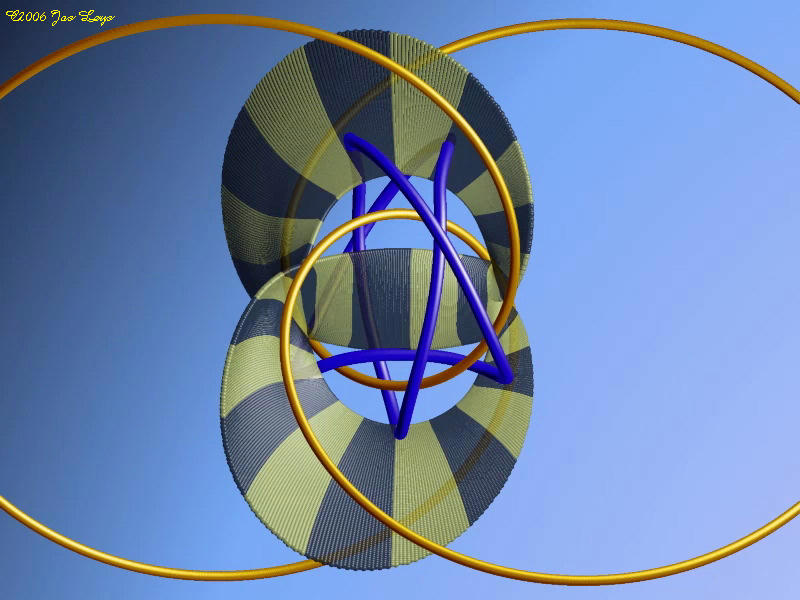}
    \includegraphics[width=0.42\textwidth]{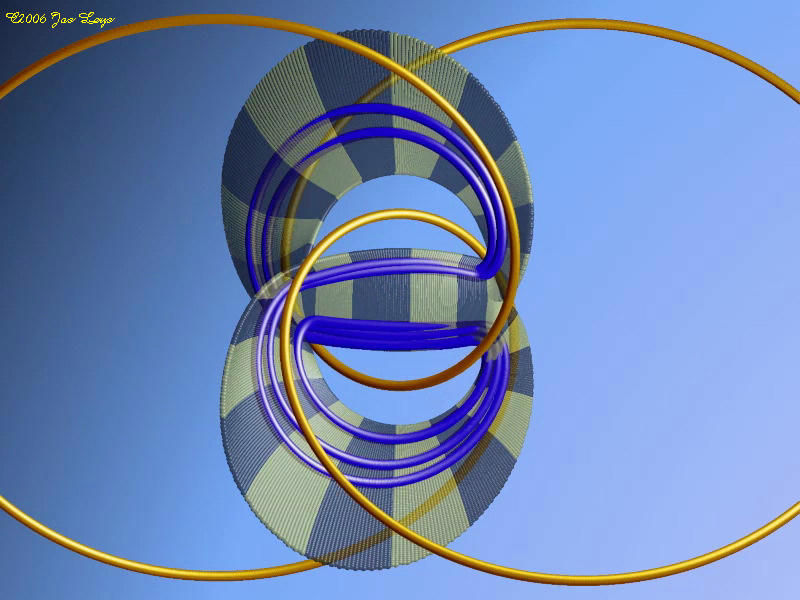}
    \caption*{Isotopy of a modular knot to a Lorenz knot. The trefoil is in yellow.}
    \label{fig:Modular-knot-in-Lorenz-Template}
\end{figure}

\'E. Ghys concludes his paper \cite{Ghys_knots-dynamics_2006} by asking for an arithmetic interpretation of the linking pairing between modular knots.
Note that the embedded Lorenz template provides a framing for the Lorenz knots, which enables to define their \emph{self-linking number} as the linking number between two parallel copies of the knot in the Lorenz template.

\section{Linking numbers of modular knots}
\label{sec:linking-numbers}

\subsection{Invariants on pairs of conjugacy classes}
\label{subsec:F(A,B)}

The action of $\PSL_2(\Z)$ on $\H\P$ and $\Tree$ enabled us to define conjugacy invariants for pairs $(A,B)$ by comparing their stable subsets in $\H\P$ or $\Tree$.
We now explain how to average those in order to obtain functions of pairs of conjugacy classes.

\subsubsection*{Summing over double cosets}

Consider a group $\Pi$ acting on a space $\Sigma$ and a function $f$ defined on $\Sigma \times \Sigma$ with values in a commutative group $\Lambda$ which is invariant under the diagonal action of $\Pi$:
\begin{equation*}
    f\colon \Sigma \times \Sigma \to \Lambda
    \qquad 
    \forall W \in \Pi, \;
    \forall a,b \in \Sigma 
    \: \colon \:
    f(a,b)=f(W\cdot a, W\cdot b).
\end{equation*}
We define an invariant $F$ for pairs of $\Pi$-orbits $[a],[b]$ by summing $f$ over all pairs of representatives of the orbits considered modulo the diagonal action of $\Pi$.
The pairs of representatives for the orbits are parametrized by the $(U\cdot a,V\cdot b)$ for $(U,V)\in \Pi / (\Stab a) \times \Pi / (\Stab b)$, and the quotient of this set by the diagonal action of $\Pi$ by left translations is denoted $\Pi / (\Stab a) \times_\Pi \Pi / (\Stab b)$.
Consequently, the sum indexed by $(U,V) \in \left(\Pi/ \Stab a \right) \times_\Pi \left(\Pi/ \Stab b \right)$ defines our desired invariant:
\begin{equation*}
    F([a],[b]) = \sum f(U\cdot A, V\cdot b).
\end{equation*}
This can also be written as the sum over double cosets $W\in (\Stab a) \backslash \Pi / (\Stab b)$:
\begin{equation*}
    F([a],[b]) = \sum f(a,W\cdot b)
\end{equation*}
because the map $(\Pi / \Stab a) \times (\Pi / \Stab b) \to (\Stab a) \backslash \Pi / (\Stab b)$ sending $(U,V)$ to $W=U^{-1}V$ is surjective, and its fibers are the orbits under the diagonal action of $\Pi$ by left translations.

We will apply this discussion to the action of $\PSL_2(\Z)$ on itself by conjugacy to obtain invariants for pairs of primitive hyperbolic conjugacy classes.
Note that in $\PSL_2(\Z)$, the centraliser of a hyperbolic $A$ is the infinite cyclic subgroup generated by its primitive root, namely the unique primitive element with a positive power equal to $A$.
Our functions $f(a,b)$ will be expressed in terms of geometrical invariants such as $\bir(A,B)$ or $\cos(A,B)$, as well as combinatorial invariants such as $\cross(A,B)$ and $\cosign(A,B)$.
To ensure that the sum is well defined, it must have finite support or converge in a completion of $\Lambda$ for an appropriate norm, and that depends on the behaviour of $f$.

\subsubsection*{Summing over $\{L,R\}$-words}

Consider a function $f$ over the pairs of coprime primitive hyperbolic $A,B\in \PSL_2(\Z)$, which is invariant under the diagonal action of $\PSL_2(\Z)$ on itself by left conjugacy.

In order to compute the sum defining $F([A],[B])$, we may group the terms $f(UAU^{-1},VBV^{-1})$ according to the $\cosign(UAU^{-1},VBV^{-1}) \in \{-1,0,1\}$ to obtain:
\begin{equation*}
    F = F_- + F_0 + F_+
\end{equation*}
The sum $F_+$ has finite support, contained in the set of pairs of Euclidean representatives for the conjugacy classes of $A,B$.
Similarly the sum $F_-$ has finite support, which we may also index by those pairs of Euclidean representatives using the fact that $\cosign(A,B)=-\cosign(A,SBS^{-1})$. Thus for $A,B\in \PSL_2(\N)$ we have the following computable expressions:
\begin{equation*}
    F_+([A],[B])=
    \sum f\left(\sigma^iA,\sigma^jB\right)
    \qquad
    F_-([A],[B]) =
    \sum
    f\left(\sigma^iA,S(\sigma^jB)S^{-1}\right)
\end{equation*}
where the indices $i\in [1,\len A]$, $j\in [1,\len B]$ are such that $\sigma^iA$ and $\sigma^jB$ end with different letters.
One may similarly split the sum $F_0$ in two parts according to the relative orientations of the axes (interchanged by the action of $S$ on one of the components of $f$), but their index sets are infinite.

Suppose that $\cosign(A,B)=0\implies f(A,B)=0$ and $f(A,B^{-1})=\epsilon f(A,B)$ with $\epsilon\in \{\pm 1\}$. This holds for $\cross$ or $\cosign$ with $\epsilon =-1$, and for their product or their absolute values with $\epsilon =1$. Then $F_0=0$ and $F_-(A,B)=\epsilon \cdot F_+(A,{}^t\!B)$ thus $F(A,B)= F_+(A,B)+\epsilon \cdot F_+(A, {}^t\!B)$.

\subsection{Linking numbers from the action on
$(\Tree,\cord)$}

The projection of the Lorenz template yields a diagram for the Lorenz link in which all crossings are positive, and those can be enumerated using the $\{L,R\}$-cycles of the corresponding modular knots. This yields the algorithmic formula \cite[4.27]{CLS_phdthesis_2022} for computing linking numbers, which was used by Pierre Dehornoy in \cite{Dehornoy_noeuds-lorenz_2011}.

We recast it in Proposition \ref{Prop:algo-sum} in terms of the action of $\PSL_2(\Z)$ on $(\Tree,\cord)$, after introducing the appropriate quantity to be summed.

\begin{Definition}
For oriented bi-infinite combinatorial geodesics $a,b \subset \Tree$ with distinct ends we define:
\begin{equation*}
    \crocs(a, b)
    =\left(\across \times \frac{1+\cosign}{2}\right)(a, b)
    =\left(\frac{1+\
    \cross}{2} \times \frac{1+\cosign}{2}\right)(a, b)
\end{equation*}
Hence $\crocs(a, b)=1$ only when the axes cross and their orientations coincide along the intersection.
\end{Definition}

We say that $A,B\in \PSL_2(\Z)$ are \emph{coprime} when their positive powers are never conjugate.

\begin{Proposition}[Algorithmic formula: sum over $\{L,R\}$-words]
\label{Prop:algo-sum}
For coprime hyperbolic elements $A,B\in \PSL_2(\N)$ we have:
\begin{equation*}
    \lk(A,B)=\tfrac{1}{2}\sum \crocs(\sigma^iA, \sigma^jB)
\end{equation*}
ranging over all $i\in [1,\len A]$, $j\in [1,\len B]$ are such that $\sigma^iA$ and $\sigma^jB$ end with different letters.
\end{Proposition}

\begin{proof}
The monoid $\PSL_2(\N)$ is endowed with the lexicographic order extending $L<R$.
The crossings between the Lorenz knots associated to $A,B\in \PSL_2(\N)$ are in bijection with the pairs of Euclidean representatives whose last letters are in the opposite order of the words themselves, so either $\sigma^iA=w_AL$, $\sigma^jB=w_BR$ with $w_A>w_B$ or $\sigma^iA=w_AR$, $\sigma^jB=w_BL$ with $w_A<w_B$.
\end{proof}

We deduce from the previous paragraph a group theoretical formula in terms of double cosets.

\begin{Theorem}[Algebraic formula: sum over double cosets]
\label{Prop:algebra-sum}
For coprime primitive hyperbolic elements $A,B\in \PSL_2(\Z)$:
\begin{equation*}
    \lk(A,B)=\tfrac{1}{2}\sum \crocs(\tilde{A}, \tilde{B})
\end{equation*}
where the sum extends over pairs of representatives $\tilde{A}=UAU^{-1}$ and $\tilde{B}=VBV^{-1}$ for the conjugacy classes with 
$(U,V)\in \Gamma/\langle A\rangle \times_\Gamma \Gamma/\langle B\rangle$. 
\end{Theorem}

\begin{Remark}[Intersection]
\label{Rem:intersection-from-link}
We recover the intersection number of modular geodesics as:
\begin{equation*}
    \lk(A,B)+\lk(A,B^{-1})=\tfrac{1}{2}\sum \across(\tilde{A},\tilde{B}) = \tfrac{1}{2}\cdot I(A,B)
\end{equation*}
whereas the sum of the cosign over pairs of intersecting axes yields:
\begin{equation*}
    \lk(A,B)-\lk(A,B^{-1})=\tfrac{1}{2}\sum \left(\across \times \cosign\right)(\tilde{A},\tilde{B}).
\end{equation*}
We deduce an efficient algorithm computing the intersection number $I(A,B)$ from $\{L,R\}$-factorisations of $A,B$ by applying algorithmic formula to the linking numbers $\lk(A,B)$ and $\lk(A,B^{-1})$.

Note that if $A$ is conjugate to $B$, then $I(A,B)$ is the intersection number between two parallel copies of the modular geodesic, which is twice its self-intersection number (counted as the number of double points).
For instance, the modular geodesic corresponding to $RLL$ has self-intersection \[\tfrac{1}{2}I([RLL],[RLL])=\lk([RLL],[RLL])+\lk([RLL],[LLR])=\tfrac{1}{2}4+\tfrac{1}{2}2=3.\]
\end{Remark}

\section{Linking function on the character variety and its boundary}

\subsection{Linking function on the character variety and its boundary}

Recall the family of Fuchsian representations $\Bar{\rho}_q\colon \PSL_2(\Z) \to \PSL_2(\R)$ parametrized algebraically by $q\in \R^*$.

\begin{Definition}
For primitive hyperbolic $A,B\in \PSL_2(\Z)$, we define the algebraic functions of $q$:
\begin{equation}\label{eq:Link}\tag{$\Link_q$}
    \Link_q(A,B)
    = \tfrac{1}{2} \sum \left(\tfrac{\asrt{\bir>1}}{\bir}\right)\left(\tilde{A}_q, \tilde{B}_q\right)
\end{equation}
\begin{equation}\label{eq:Cos}\tag{$\Cos_q$}
    \Cos_q(A,B)
    = \tfrac{1}{2} \sum (\across \times \cos)(\tilde{A}_q, \tilde{B}_q)
\end{equation}
by summing over the pairs of representatives $\tilde{A}=UAU^{-1}$ and $\tilde{B}=VBV^{-1}$ for the conjugacy classes of $A$ and $B$ where $(U,V)\in \Gamma/\Stab(A) \times_\Gamma \Gamma/\Stab(B)$.
\end{Definition}

The appearance of $\asrt{\bir>1}=\across$ as a factor in the terms of \ref{eq:Link} and \ref{eq:Cos} amounts to restricting the summations over pairs of matrices whose axes intersect.
Hence the support of the sums corresponds to the intersection points of the modular geodesics $[\alpha]$ and $[\beta]$ associated to the conjugacy classes, which must be counted with appropriate multiplicity when $A$ or $B$ is not primitive. Thus:
\begin{equation*}
    \Link_q(A,B)
    = \tfrac{1}{2} \sum \left(\cos \tfrac{\theta}{2}\right)^2 
    \qquad \mathrm{and} \qquad
    \Cos_q(A,B)
    = \tfrac{1}{2} \sum \left(\cos \theta\right).
\end{equation*}
Observe that since $\bir(A,B)^{-1}+\bir(A,B^{-1})^{-1}=1$ we have $\Link_q(A,B)+\Link_q(A,B^{-1})=\tfrac{1}{2}I(A,B)$.

\begin{Conjecture}
The angles turning from $[\alpha_q]$ to $[\beta_q]$ in the direction prescribed by the orientation of $\M_q$ have cosines $(\cross \times \cos)(\tilde{A_q},\tilde{B_q})$: we believe that they sum up to $0$, as explained in \ref{subsec:Link-Fuchsian-group}.
\end{Conjecture}

\begin{Theorem}
\label{Thm:Bir(A,B)-->lk(A,B)}
For primitive hyperbolic conjugacy classes $[A],[B]$ in $\PSL_2(\Z)$ we have:
\begin{align*}
    &\Link_q(A,B) \xrightarrow[q\to \infty]{} \lk(A,B)
    \\
    &\Cos_q(A,B) \xrightarrow[q\to \infty]{} 
    \lk(A,B)-\lk(A^{-1},B)
    =2\lk(A,B)-\tfrac{1}{2}I(A, B)
\end{align*}
\end{Theorem}

\begin{proof}
Recall from Lemma \ref{Lem:cos-cosh-sinh} the relation $1/\bir(A_q,B_q)=\tfrac{1}{2}(1+\cos(A_q,B_q))$ and from Lemma \ref{Lem:cos_q-limit} the limit $\cos(A_q,B_q)\to \cosign(A,B)$ as $q\to \infty$. Hence the terms of the sum defining \ref{eq:Link} converge to those in the sum of Proposition \ref{Prop:algebra-sum}. The limit of \ref{eq:Cos} follows from Remark \ref{Rem:intersection-from-link}.
\end{proof}

Let us display the graphs of \textcolor{blue}{$2\Link_q(A,B)$} and \textcolor{red}{$2\Link_q(A,B^{-1})$} along with their average \textcolor{black!50!green}{$\tfrac{1}{2}I(A,B)$} as a function of $q\in \R$ for some $A,B\in \PSL_2(\N)$. The legend $A=[a_0,a_1,\dots]$ means $A=R^{a_0}L^{n_1}\cdots$.

\begin{figure}[h]
    \centering
    \includegraphics[width=0.36\textwidth]{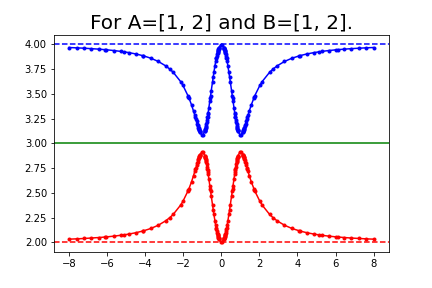}
    \hspace{-0.9cm}
    \includegraphics[width=0.36\textwidth]{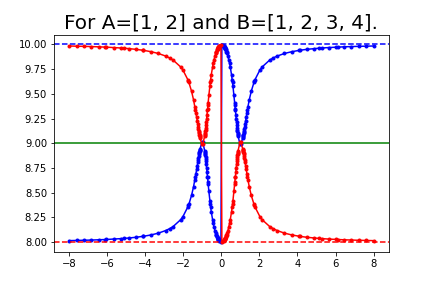}
    \hspace{-0.9cm}
    \includegraphics[width=0.36\textwidth]{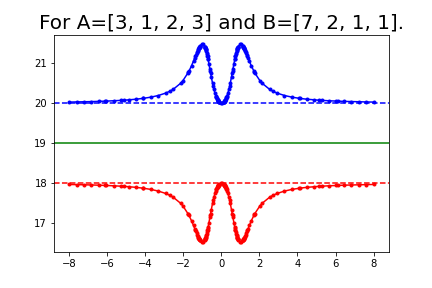}
    \vspace{-0.4cm}
    \caption*{\textcolor{blue}{$L_q(A,B)$} interpolates between the arithmetic at $q=1$ and the topology at $q=\infty$.}
\end{figure}

\subsection{Graphs of $\Link_q(A,B)$ for $q\in \C$}

The function $q\mapsto \Link_q(A,B)$ is the ratio of a polynomial by the square-root of an even polynomial $\sqrt{\disc(A_q)\disc(B_q)}$, so its extension to $q\in \C$ is fixed by choosing the standard branch of the square-root.
We represent the complex graph of $q\mapsto \Link_q(A,B)$ by assigning a colour to each point of $\C$ using the HSV colour scheme: the hue varies with the argument and the brightness varies with the modulus.
Since $\Link_q=\Link_{1/q}$ we restrict $\lvert q \rvert < 1+\epsilon$ for some $\epsilon>0$ chosen according to aesthetic criteria.

\begin{figure}[h]
    \centering
    \vspace{-0.2cm}
    \hspace{-2cm}
    \includegraphics[width=0.5\textwidth]{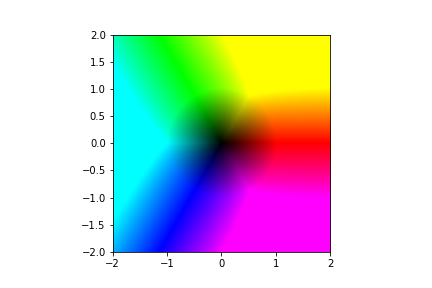}
    \hspace{-3.4cm}
    \includegraphics[width=0.5\textwidth]{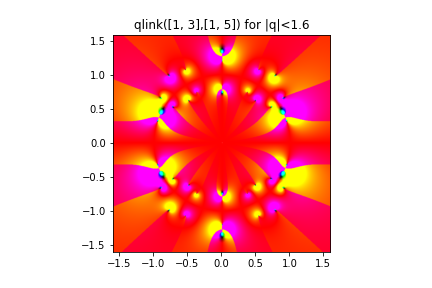}
    \hspace{-3.4cm}
    \includegraphics[width=0.5\textwidth]{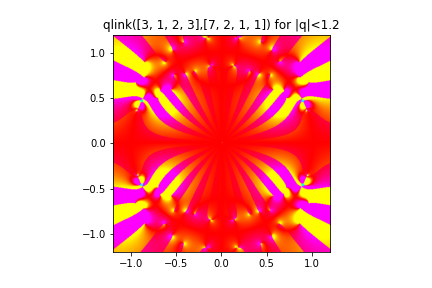}
    \hspace{-3.4cm}
    \\
    \hspace{-2cm}
    \includegraphics[width=0.5\textwidth]{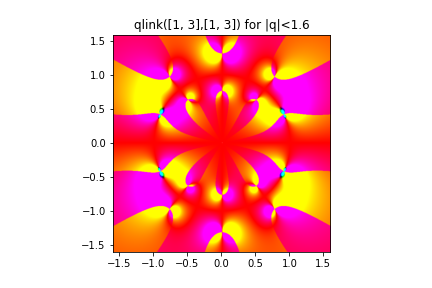}
    \hspace{-3.4cm}
    \includegraphics[width=0.5\textwidth]{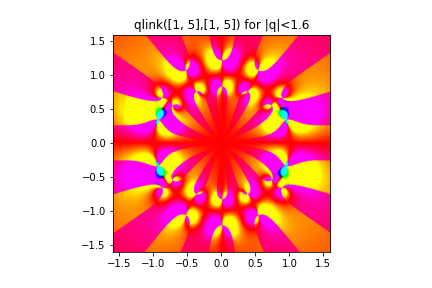}
    \hspace{-3.4cm}
    \includegraphics[width=0.5\textwidth]{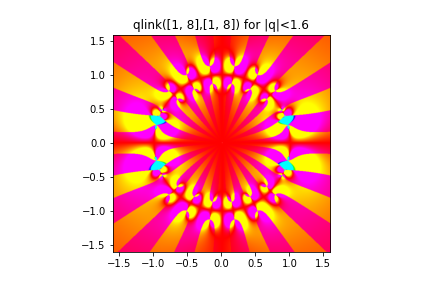}
    \hspace{-3.4cm}
    \\
    \hspace{-2cm}
    \includegraphics[width=0.5\textwidth]{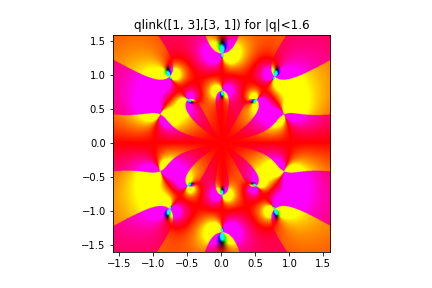}
    \hspace{-3.4cm}
    \includegraphics[width=0.5\textwidth]{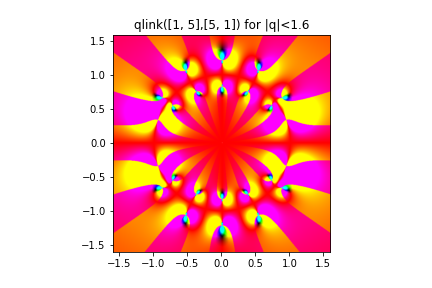}
    \hspace{-3.4cm}
    \includegraphics[width=0.5\textwidth]{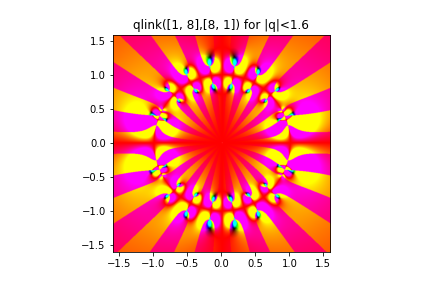}
    \hspace{-3.4cm}
    \\
    \vspace{-0.4cm}
    \caption*{The identity map and some graphs of $\Link_q(A,B)$ for $q\in \C$ with $\lvert q\rvert < 1+\epsilon$.}
\end{figure}

The main observation is that the zeros and poles of $\Link_q(A,B)$ seem to concentrate on the unit circle.
This is neither surprising nor obvious, as we explain in the next paragraph. 

\subsection{Locating poles and zeros of \texorpdfstring{$\Link_q$}{L_q}} 
\label{subsec:L_q(A,B)=0}

Let us first relate \cite[Proposition 5.16]{CLS_phdthesis_2022}.

A primitive hyperbolic conjugacy class in $\PSL_2(\Z)=\pi_1(\M)$ corresponds to a primitive modular geodesic in $[\alpha]\subset \M$. It lifts to a modular knot in $[\Vec{\alpha}] \subset \U$ which in turn yields a conjugacy class in $\BB_3=\pi_1(\U)$.
A conjugacy class in the braid group on three strands defines, by taking its closure, a link $\sigma_A$ in a solid torus.
In \cite[Proposition 5.16]{CLS_phdthesis_2022} we relate the Alexander polynomial $\Delta(\sigma_A)\in \Z[t^{\pm 1}]$ of this link $\sigma_A$ to the Fricke polynomial $\Tr A_q \in \Z[q^{\pm 1}]$ of the modular geodesic $[\alpha]$.

\begin{Proposition}[Alexander, Fricke and Rademacher]
\label{Prop:Fricke-Alex}
For a primitive hyperbolic $A\in \SL_2(\N)$, the Alexander polynomial of the link $\sigma_A$ is given in terms of $q=\sqrt{-t}$ by:
\begin{equation*}
    \Delta(\sigma_A)=\left(q^{\Rad(A)}-\Tr(A_q)+q^{-\Rad(A)}\right) \div (q-q^{-1})^2
\end{equation*}
\end{Proposition}

\begin{proof}
    We sketch the ideas of the proof and refer to \cite[Section 4.1 and Section 5.2]{CLS_phdthesis_2022} for details.
    
    We first show (in \cite[Theorem 4.16]{CLS_phdthesis_2022}, relying on the description of the Lorentz template), that the geometric section of loops up the Seifert fibration $\U\to \M$ defined by lifting modular geodesics corresponds to the combinatorial section of conjugacy classes up the central extension $\BB_3\to \PSL_2(\Z)$ defined by the morphism of monoids $(L,R)\mapsto (\sigma_1^{-1}, \sigma_2)$.
    In other terms, the modular geodesic $[\alpha]\subset \M$ associated to the conjugacy class in $\PSL_2(\Z)$ encoded by the $\{L,R\}$-cycle $A\in \PSL_2(\N)$ lifts to the modular knot $[\vec{\alpha}]\subset \U$ yielding a conjugacy class $\BB_3$ associated to the same $(\sigma_1^{-1}, \sigma_2)$-cycle $\sigma_A$.
    Besides, following \cite{Birman-Brendle_survey-braids_2005}, the Alexander polynomial $\Delta(\sigma)$ of a braid $\sigma\in \BB_3$ is given in terms of the Burau representation $\operatorname{Br}_t\colon \BB_3 \to \GL_2(\Z[t^{\pm 1}])$ defined by: 
    \begin{equation*}
    \operatorname{Br}_t(\sigma_1)=\begin{psmallmatrix}
        -t & 1 \\ 1 & 0
    \end{psmallmatrix}
    \quad \mathrm{and} \quad
    \operatorname{Br}_t(\sigma_2)=\begin{psmallmatrix}
        1 & 0 \\ t & -t
    \end{psmallmatrix}
    \qquad \mathrm{as}\qquad
    \Delta(\sigma)(t)=\tfrac{\det(\operatorname{Br}_t(\sigma))}{1+t+t^2}.
    \end{equation*}
    Conjugating $\operatorname{Br}_t$ by $S$ and setting $q=\sqrt{-t}$ yields the representation $\operatorname{Sq}\colon \BB_3\to \SL_2(\Z[q^{\pm1}])$ given by $\sigma_1^{-1} \mapsto \tfrac{1}{q}L_q$ and $\sigma_2 \mapsto q R_q$. 
    Thus for $A\in \SL_2(\N)$ we have $\operatorname{Sq}(\sigma_A)=q^{\Rad(A)}\rho_q(A)$. 
    We conclude using the Cayley-Hamilton relation $\det(M-\Id)=\det(M)-\Tr(M)+1$ for $M\in \Gl_2$, and replacing $1+t+t^2=q(q-q^{-1})^{2}$, noting that $\Delta(\sigma)$ is defined up to multiplication by a unit $t^{\pm n}\in \Z[t^{\pm 1}]$.
\end{proof}

\begin{Question}[The poles of \texorpdfstring{$\Link_q$}{L_q}]
Locating the zeros of Alexander polynomials of various classes of knots and links has been the subject of various conjectures and results.
For instance \cite{Stoimenov_alexander-hoste-conjecture_2019} studies a conjecture of Hoste for links with braid index $3$.
Now recall that $\Link_q(A,B)-\Link_q(A,B^{-1})$ can be expressed as a finite sum of terms 
of the form $\cos(A_q,B_q)$.
This is why one may guess with Proposition \ref{Prop:Fricke-Alex} a concentration property for the zeros and poles of $\Link_q$ around the unit circle.
Still, it would remain a challenge to prove it.
\end{Question}

\begin{Remark}
We should also mention that \cite{Dehornoy_zeros-alex-modular-knots_2015} has shown concentration properties for the zeroes of the Alexander polynomials of Lorenz knots: they lie an annulus whose inner and outer radii are bounded in terms the genus and the braid index of the knot. 
(Beware not to confuse the Alexander polynomials of modular knots $[\Vec{\alpha}]$ and the Alexander polynomials of modular braids $\sigma_A$.)

We may wonder about the relation between the the links consisting of the trefoil with a modular knot and the $3$-braids.
In a similar line of thought, one may ask what the linking polynomial $\Link_q(A,B)$ says about the modular link consisting of $[\Vec{\alpha}]$, $[\Vec{\beta}]$ and the trefoil.

Let us mention how modular links and $3$-braids fit into one picture.
Consider the tautological bundle of unit covolume lattices over the base $\SL_2(\Z)\backslash \SL_2(\R)$, whose points are pairs consisting of a lattice $\Lambda \subset \R^2$ and a matrix $M\in \SL_2(\R)$ such that $\Lambda=\Z^2 M$.
We can define the Weierstrass points in every fiber by considering the lattice $\tfrac{1}{2}\Lambda \supset \Lambda$, and this yields a "Weierstrass $4$-valued section".
While going around a modular knot $[\Vec{\alpha}]$ downstairs, the monodromy along this Weierstrass section in the fiber corresponds to (a double cover) of the braid $\sigma_A$.
\end{Remark}

\section{Linking numbers and homogeneous quasimorphisms}

\label{sec:quasi-morphism}

\subsection{Combinatorial formula: sum of linked patterns}

We now derive a combinatorial formula for the linking numbers \cite[Proposition 4.34]{CLS_phdthesis_2022} arising from a different count of the crossings in the Lorenz template. It follows from the algorithmic formula, but we propose a visual proof.

The monoid $\PSL_2(\N)$ freely generated by $\{L,R\}$ has the lexicographic order extending $L<R$. 
The monoid $\PSL_2(\N)\setminus\{\Id\}$ maps to the set $\{L,R\}^{\N}$ of infinite binary sequences by sending a finite word $A$ to its periodisation $A^\infty$. This map is increasing and injective in restriction to primitive words.
We denote by $\sigma$ the Bernoulli shift on $\{L,R\}^\N$ which removes the first letter, as well as the cyclic shift on $\PSL_2(\N)\setminus\{\Id\}$ which moves the first letter at the end.
These shifts are intertwined by the periodisation map: $(\sigma^j A)^\infty = \sigma^j(A^\infty)$.
For a pattern $P\in \PSL_2(\N)$ and an infinite order $A\in \PSL_2(\N)$, let $\pref_P(A^\infty)=\asrt{A^\infty\in P\cdot\PSL_2(\N)}\in \{0,1\}$ tell whether $P$ is a prefix of $A^\infty$, and $\occ_P(A) = \sum_{j=1}^{\len A}
\pref_P\left(\sigma^jA^{\infty} \right)$ count the number of cyclic occurrences of $P$ in $A\bmod{\sigma}$.

Recall that $A,B\in \PSL_2(\Z)$ are coprime when their positive powers are never conjugate.
Thus $A,B\in \PSL_2(\N)$ are not coprime when they admit cyclic permutations generating submonoids with non-trivial intersection, in other terms if $A^\infty = B^\infty \bmod{\sigma}$.

\begin{Proposition}[Combinatorial formula: sum of linked patterns]
\label{Prop:sum-linked-patterns}
For coprime hyperbolic $A,B\in \PSL_2(\N)$ the corresponding modular knots have linking number:
\begin{equation}
\label{eq:sum-linked-patterns}\tag{SLP}
\lk(A,B) = \frac{1}{2} \sum_{w}
\begin{pmatrix}
\occ_{RwL}(A)\cdot \occ_{LwR}(B)
\\+\\
\occ_{LwR}(A)\cdot \occ_{RwL}(B)
\end{pmatrix}
\end{equation}
where the summation extends over all words $w\in \PSL_2(\N)$ including the empty one.
\end{Proposition}

\begin{proof}[Visual proof \texorpdfstring{\cite[Figure 4.9]{CLS_phdthesis_2022}}{VisualProof49}]
Split the Lorenz template by extending the dividing line backwards in time, and observe the crossings appearing in its standard planar projection: they occur in regions arranged according to a binary tree indexed by pairs of words of the form $(RwL,LwR)$.
\end{proof}

\begin{Remark}
\label{rem:long-patterns}
For $\len(P)\ge \len(A)$ we have $\occ_P(A)>0$ if and only if $A^\infty = P^\infty \bmod{\sigma}$, which is equivalent to the non-coprimality of $P$ and $A$.
Hence the coprimality assumption on $A$ and $B$ ensures that the support of the sum \eqref{eq:sum-linked-patterns} is contained in the set of $w$ such that $\len w < \max\{\len A,\len B\}$.

If $A$ and $B$ are not coprime, then they are conjugate to positive powers $C^m$ and $C^n$ of a primitive $C\in \PSL_2(\N)$, and restricting the sum \eqref{eq:sum-linked-patterns} to the indices $w$ with $\len w < \max\{\len A,\len B\}$ yields $mn$ times the self-linking number of the modular knot associated to $C$ with the Lorenz framing.
\end{Remark}

\begin{Remark}
The cycle $A \bmod{\sigma}$ has a multiset of $L$-exponents and a multiset of $R$-exponents. 
Formula \eqref{eq:sum-linked-patterns} shows that $\lk(A,RL^{m+1})-\lk(A,RL^{m})$ counts the number of $L$-exponents which are $>m\ge 1$ and that $\lk(A,LR^{n+1})-\lk(A,LR^{n})$ counts the number of $R$-exponents which are $>n\ge 1$.
\end{Remark}

\begin{Theorem}[Link equivalence]
\label{Thm:linkeq_implies_conjugate}
If hyperbolic $A,B\in \PSL_2(\Z)$ are link equivalent, namely $\lk(A,C)=\lk(B,C)$ for all hyperbolic $C \in \PSL_2(\Z)$, 
then they are conjugate.
\end{Theorem}

\begin{proof}

The set $\mathcal{Z}=R\PSL_2(\N)L\sqcup L\PSL_2(\N)R$ of $\{L,R\}$-words which start and end with different letters is endowed with the involution $z\mapsto \Bar{z}$ exchanging the extreme letters, having no fixed points.
The free $\Z$-module $\Omega$ generated by the set $\mathcal{Z}$ is naturally decomposed as the direct sum of rank-$2$ sub-modules generated by pairs $\{z,\Bar{z}\}$.
It is therefore endowed with a non-degenerate symmetric bilinear form $\Omega \times \Omega \to \Z$, given by the direct sum of the hyperbolic structures on those planes:
\begin{equation*}
    \Omega 
    = \bigoplus_{z\in \mathcal{Z}} \Z\cdot z
    = \bigoplus_{z>\Bar{z}} \Z\cdot z \oplus \Z\cdot \Bar{z}
    \qquad
    (a\cdot b) 
    = \sum_{z\in \mathcal{Z}} a_z b_{\Bar{z}} 
    = \sum_{z>\Bar{z}} (a_z b_{\Bar{z}} + a_{\Bar{z}} b_{z})
\end{equation*}

The length function $\len \colon \PSL_2(\N) \to \N$ yields a filtration of the set $\mathcal{Z}$ by the chain of subsets $\mathcal{Z}_n = \{z \in \mathcal{Z} \mid \len(z)\le n\}$ with cardinals $2^{n-1}$, which is invariant by the involution.
This induces a filtration of the module $\Omega$ by the corresponding chain of sub-modules $\Omega_n$ with ranks $2^{n-1}$, all invariant under the orthogonal symmetry.
Thus each unimodular quadratic $\Z$-module $\Omega_n$ (decomposed as a direct sum of hyperbolic planes) is canonically isomorphic to its dual $\Omega_n^*$.

Now consider cyclic words $A,B\in \PSL_2(\N) \bmod{\sigma}$ corresponding to link equivalent hyperbolic conjugacy classes in $\PSL_2(\Z)$. 
Let $m=\max\{\len(A),\len(B)\}$ and consider the linear forms on $\Omega_{m}$ defined by the sequences $(\occ_z(A))_{z}$ and $(\occ_z(B))_z$ for $z\in \mathcal{Z}_m$.
Since $A,B$ are link equivalent, the isomorphism $\Omega_m\to \Omega_m^*$ implies by Theorem \ref{Thm:linkeq_implies_conjugate} and Remark \ref{rem:long-patterns} that these linear forms coincide, so that $\occ_z(A)=\occ_z(B)$ for all $P\in \mathcal{Z}_m$.
In particular, for $z$ a linear representative of $B$ we find that $\occ_B(A)=\occ_B(B)>0$ whereby $A=B \bmod{\sigma}$.
\end{proof}

\subsection{Homogeneous quasi-morphisms on the modular group}

For a group $\Pi$, a function $f\colon \Pi \to \R$ is called a \emph{quasi-morphism} if it has a bounded derivative:
\[df(A,B)=f(B)-f(AB)+f(A).\]
A quasi-morphism $f\colon \Pi \to \R$ is called \emph{homogeneous} if it is a morphism in restriction to the abelian subgroups of $\Pi$ (which for $\Pi = \PSL_2(\Z)$ means that $f(A^n)=nf(A)$ for all $A\in \Pi$ and $n\in \N$).
Observe that a homogeneous quasi-morphism is bounded only if it is trivial, necessarily constant on conjugacy classes, and vanishes on torsion classes.
The real vector space $PX(\Pi)$ of homogeneous quasi-morphisms is a Banach space for the norm $\lVert df \rVert_\infty$, as was shown in \cite{MatsuMorita_Hb(Homeo)_1985, Ivanov_H2b(G)-Banach_1988}.

For a pattern $P\in \PSL_2(\N)$ we define the $P$-asymmetry of an infinite order $A\in\PSL_2(\N)$ by \[\mas_P(A)=\occ_P(A)-\occ_{{}^t\!P}(A)\]
Notice that $\mas_P(A)=\occ_P(A)-\occ_P({}^t\!A)$ and that ${}^t\!A$ is conjugate to $A^{-1}$ by $S\in \PSL_2(\Z)$.
Extending $\mas_P(A)=0$ for elliptic $A$ yields a conjugacy invariant function $\mas_P \colon \PSL_2(\Z)\to \Z$.
In particular for $P=R$ we recover the \emph{Rademacher function} as $\mas_P(A) = \Rad(A)$.

\begin{Lemma}
\label{Lem:cocycle}
For all $P\in \PSL_2(\N)$, the function $\mas_P \colon \PSL_2(\Z)\to \Z$ is a homogeneous quasi-morphism. If $P\ne {}^t\!P$ then $\mas_P$ is unbounded, and if $P$ does not overlap itself then $\lVert d\mas_P \rVert_\infty \le 6$.
\end{Lemma}

\begin{proof}
The proof relies on the ideas in \cite{BarGhys_cocycles-actions-arbres_1991} (see also \cite[Lemma 5.3]{Grigorchuk_bounded-cohomology_1995}).
\end{proof}

\begin{Theorem}
\label{Thm:Cos_A}
For every hyperbolic $A\in \PSL_2(\Z)$, the function $\Cos_A\colon B\mapsto \lk(A,B)-\lk(A^{-1}, B)$ is a homogeneous quasi-morphism $\PSL_2(\Z)\to \Z$, which is unbounded unless $A$ is conjugate to $A^{-1}$.
It can be computed for $A,B\in \PSL_2(\N)$ as:
\begin{equation}
\label{eq:Cos_A} \tag{$\Cos_A$}
\Cos_A(B) = \lk(A,B)-\lk(A, {}^t\!B) 
= \frac{1}{2} \sum_{w}
\begin{pmatrix}
\occ_{RwL}(A)\cdot \mas_{LwR}(B)
\\+\\
\occ_{LwR}(A)\cdot \mas_{RwL}(B) 
\end{pmatrix}
\end{equation}
where the summation extends over all words $w\in \PSL_2(\N)$ with $\len(w)<\max\{\len A, \len B\}$.
\end{Theorem}

\begin{proof}
The quantity $\Cos_A(B)$ is homogeneous in $A$ and $B$.
Let us explain why $\Cos_A$ is a quasi-morphism for $A$ primitive.
Recall that $A^{-1}$ and ${}^t\!A$ are conjugate by $S$ and notice that $\lk(A^{-1},B)=\lk(A, B^{-1})$.
Therefore \eqref{eq:sum-linked-patterns} yields \eqref{eq:Cos_A} and $d\Cos_A=\tfrac{1}{2}\sum_w \left(\occ_{RwL}(A) \cdot d\mas_{LwR}+\occ_{LwR}(A) \cdot d\mas_{RwL}\right)$.

Since $\occ_P(A)\le \len(A)$ it is enough to prove by Proposition \ref{Lem:cocycle} that for every $X,Y\in \PSL_2(\Z)$, the sum $d\Cos_A(X,Y)$ contains at most $\len(A)^2$ non-zero terms with $\len(w)\ge \len(A)$.

The $\{L,R\}$-words $w$ with $\occ_{LwR}(A)>0$ correspond to the triples $1\le m,n\le \len(A)$ and $k\in \N$ such that $\sigma^mA=LuRv$, $\sigma^nA=RvLu$, and $LwR=L(uRvL)^kuR$ for some $\{L,R\}$-words $u,v$. In this situation we write $P_{mn}^k=LwR=L(uRvL)^kuR$ and $Q_{mn}^k=RwL=R(uRvL)^kuL$.

By construction (and the primitivity of $A$), two distinct $Q_{mn}^k$ cannot overlap except along a prefix and suffix of length $<\len(A)$.
We may thus adapt the argument for \cite[Proposition 5.10]{Grigorchuk_bounded-cohomology_1995}.
The quantity $d\mas_Q(X,Y)$ measures the "$Q$-perimeter" of a tripod $(*,X*,XY*)$ in the tree $\Tree$, and it is non-zero only if $Q$ can be matched along a portion covering its incenter.
But for each $(m,n)$, at most two values of $k>1$ may lead to such patterns $Q_{mn}^k$.
The same reasoning applies with $L, R$ interchanged.
This proves the bound on the number of non-zero summands for $d\Cos_A$.

Finally by Theorem \ref{Thm:linkeq_implies_conjugate} we have $d\Cos_A=0$ only if $A$ is conjugate to ${}^t\!A$.
\end{proof}

Let $\mathcal{P}$ denote the set of \emph{Lyndon words}, namely the $\{L,R\}$-words which are greater than every one of their cyclic permutations.
Notice that such words are primitive, and cannot overlap themselves.
Hence if a Lyndon word is equal to a cyclic permutation of its transpose then it is actually symmetric. Let $\mathcal{P}_0$ be the subset of symmetric Lyndon words and choose a partition $\mathcal{P}\setminus \mathcal{P}_0=\mathcal{P}_-\sqcup \mathcal{P}_+$ in two subsets which are in bijection by the transposition.
Note that $\mathcal{P}$ indexes the set of primitive infinite order conjugacy classes in $\PSL_2(\Z)$, and $\mathcal{P}_0$ the subset of those which are stable under inversion.
Of course $\Id \in \mathcal{P}_0$, we may choose $R\in \mathcal{P}_+$, and denote $\Cos_R:=\mas_R=\Rad$ by convention.

\begin{Proposition}
The collection of $\Cos_A\in PX(\PSL_2(\Z);\R)$ for $A\in \mathcal{P}_+$ is linearly independent.
\end{Proposition}

\begin{proof}
Consider distinct Lyndon words $A_1,\dots,A_k\in \mathcal{P}_+\setminus\{R\}$, and let $A_1,\dots,A_j$ be those of maximal length $m$.
The $\Cos_{A_j}$ are linearly independent of $\Cos_R$ as $\Cos_{A_j}(R)=0$ whereas $\Cos_R(R)=1$.
Suppose by contradiction that we have a linear relation $\sum r_i \Cos_{A_i}=0$ for $r_i\in \R^*$.
As in the proof of Proposition \ref{Thm:linkeq_implies_conjugate}, this restricts to a linear relation in $\Omega_m^*$, and using the isomorphism $\Omega_m\to \Omega_m^*$ we find that $\sum r_i \occ_P(A_i) = \sum r_i \occ_P({}^t\!A_i)$ for all $P\in \mathcal{Z}_m$.

For $P=A_j$, we have $\occ_{P}(A_i)=0=\occ_{P}({}^t\!A_i)$ for all $i>j$ because $A_j$ cannot overlap itself, and $\occ_{P}(A_i)=0=\occ_{P}({}^t\!A_i)$ for $i<j$ because the Lyndon word $A_j$ is different from $A_i$.
We also have $\occ_{P}({}^t\!A_j)=0$ since $A_j$ is not conjugate to its inverse, whence not link equivalent to its transpose by Proposition \ref{Thm:linkeq_implies_conjugate}.
Since $\occ_{P}(A_j)=1$ we have $r_1=0$, which is the desired contradiction.
\end{proof}

\begin{Proposition}
Every $f\in PX(\PSL_2(\Z))$ can be written $\sum c_A(f) \cdot \Cos_A$ for unique $c_A(f) \in \R^{\mathcal{P}_+}$.
\end{Proposition}

\begin{proof}
A homogeneous quasi-morphism $f\in PX(\PSL_2(\Z))$ quotients to a function on the set of infinite order primitive conjugacy classes, or on $\mathcal{P}$. It must vanish on $\mathcal{P}_0$, and change sign by transposition, so the function restricted to $\mathcal{P}_+$ uniquely determines $f$.
Since $\Cos_A(R)=0$ for all $A\in \mathcal{P}_+\setminus \{R\}$ and $\Cos_R(R)=1$ we may assume $f(R)=0$ from now on.

Recall the structure of the filtered quadratic space $\Omega$ introduced in the previous proofs.
Notice that the cyclic shift acting on $\mathcal{Z}$ preserves the $\mathcal{Z}_n$ (but does not commute with the involution $z\mapsto \Bar{z}$ for $n>2$ as these two actually generate the full group of permutations $\mathfrak{S}_n$).
Fix $m\in \N$ and consider the subspace $\Lambda^*_n \subset \Omega^*_m$ of elements which are invariant under the shift $\sigma$ and change sign under transposition.
Its elements are uniquely determined by their values on $\mathcal{L}_m := \mathcal{Z}_m \cap \mathcal{P}_+$. 
It contains the $(\mas_P(z))_{z\in \mathcal{Z}}$ for $P\in \mathcal{L}_m$, as well as the $(\Cos_A(z))_{z\in \mathcal{Z}}$ for $A\in \mathcal{L}_m$.
We know from the previous proof that the latter is free, and can be expressed as linear combinations of the former, so both of these form bases of $\Lambda^*_m$, whose dimension equals the cardinal of $\mathcal{L}_m$.

Hence the restriction $f_m\in \Lambda^*_m$ of $f$ to $\mathcal{L}_m$ can be expressed as a linear combination of the $\mas_P$ or of the $\Cos_A$.
We thus have a projective system of elements $f_m$ in the linear vector spaces $\Lambda^*_{n}$ with compatible bases so the coefficients of the limit $f=\varprojlim f_m$ are well defined in either basis.
\end{proof}

In passing, we recovered the following reformulation of \cite[Theorem 5.11]{Grigorchuk_bounded-cohomology_1995}.

\begin{Corollary}
\label{Cor:mes_P-basis}
The collection of $\mas_P\in PX(\PSL_2(\Z))$ for $P\in \mathcal{P}_+$ is linearly independent, and every $f\in PX(\PSL_2(\Z);\R)$ can be written $\sum m_P(f) \mas_P$ for unique $m_P(f) \in \R^{\mathcal{P}_+}$.
\end{Corollary}

\section{Further directions of research}


\subsection{Linking forms of Fuchsian groups}
\label{subsec:Link-Fuchsian-group}

To begin with, we may compare the definitions of the functions $\Link_q$ and $\Cos_q$ and their limiting behaviour at $q=\infty$ with similar considerations which have been made for non-oriented loops in a closed surface $S$ of genus $g\ge 2$.
Such loops, corresponding to the conjugacy classes of $A,B\in \pi_1(S)$ up to inversion, define trace functions $\Tr(A), \Tr(B)$ on the $\SL_2(\C)$-character variety of $\pi_1(S)$, whose real locus contains the Teichm\"uller space of $S$ as a Zariski-dense open set.
This character variety carries a natural symplectic structure \cite{Goldman_symplectic-nature-pi1_1984}, given by the Weil-Petersson symplectic form.

The sum $\Cos_q(A,B)$ looks very much like Wolpert's cosine formula \cite{Wolpert_fenchel-nielsen-deformation_1982, Wolpert_formula-cosine-Fenchel-Nielsen_1982} computing the Poisson bracket  $\{\Tr(A),\Tr(B)\}$ of the trace functions.
In fact, Wolpert sums $\cross(A,B)\cos(A,B)$ over the intersection points $p\in [\alpha]\cap[ \beta]$, that is the cosines of the angles turning from $\alpha$ to $\beta$ in the direction prescribed by the orientation of the surface.
Hence while our cosine formula is a symmetric formula of oriented geodesics, Wolpert's cosine formula yields a skew-symmetric function of non-oriented geodesics.
Note however that the Teichm\"uller space of $\M$ is reduced to a point so any Poisson structure in the usual sense would be trivial, so in our setting we expect Wolpert's sum to be identically zero (as corroborated by our computer experimentation).

Moreover, the Weil-Petersson symplectic form has been extended to several compactifications of the character variety \cite{PapadoPenne_forme-symplectic-bord-Teichmuller_1991, Sozen-Bonahon_weil-petersson-thurston-symplectic_2001, MS_ML-Newton-Poisson_2021}.
The limits of the Poisson bracket $\{\Tr(A), \Tr(B)\}$ at the respective boundary points have been interpreted in \cite[Proposition 6]{Bonahon_earthquake-mesaured-laminations_1992} and \cite{MS_ML-Newton-Poisson_2021}. 
Thus, we may generalise the definitions of our functions \ref{eq:Link} and \ref{eq:Cos} to oriented geodesics in hyperbolic surfaces and ask for interpretations of their limits at boundary points of the Teichm\"uller space.
We believe that \ref{eq:Link} and \ref{eq:Cos} extend by continuity to pairs $A,B$ of oriented geodesic currents.
This should be analogous to the extension of the intersection form described in Bonahon \cite{Bonahon_geodesic-currents_1988}.

Pursuing this direction, one may also wish to replace $\rho$ with a representation $\Gamma \to \operatorname{Homeo}^+(\S^1)$, a metric of negative curvature, or a generalised cross-ratio \cite{Otal_symplectique-bord-birapport_1992, LabourieMcShane_cross-ratios_2009}.
The aim would be to think of \ref{eq:Link} and \ref{eq:Cos} as differential forms on the "tangent bundle" to these spaces of representations, metrics or cross-ratios, considered up to appropriate equivalence relations.
Indeed, for any group $\Pi$, the semi-conjugacy classes of representations $\Pi \to \operatorname{Homeo}^+(\S^1)$ correspond \cite{Ghys_H2b(Homeo(S1);R)_1984} to the "integral points of the unit ball" in the second bounded cohomology group $H^2_b(\Pi;\R)$, namely the elements represented by bounded $2$-cocycles with values in $\{-1,0,1\}$. For $\Pi = \pi_1(S)$ it contains \cite{BargeGhys_H2b(Surface)_1988} the space of differential $2$-forms on $S$, and we suspect that something similar is true for some spaces of generalized cross-ratios, thus we ask:

\begin{Question}
How to interpret $\Link_\rho(A,B)$ or $\Cos_\rho(A,B)$ as "differential forms" on (an appropriate subspace in) the second bounded cohomology group $H^2_b(\Gamma;\R)$ ?
\end{Question}


\subsection{Arithmetic and Geometric deformations}

Let us mention another general context in which our definitions \ref{eq:Link} and \ref{eq:Cos} seem to apply with almost no changes.
Recall that our definitions of the cross-ratio and cosine in paragraph \ref{subsec:disc-bir_K} hold for pairs of semi-simple elements in $\PGL_2(\Field)$.
Thus for any faithful representation of a group $\rho \colon \Gamma \to \PSL_2(\Field)$ sending $A,B\in \Gamma$ to semi-simple elements, one may define the following invariants for the pair of conjugacy classes: \[\Link_\rho(A,B)=\sum \bir(\rho\tilde{A},\rho \tilde{B})^{-1}\qquad \Cos_\rho(A,B)=\sum \cos(\rho\tilde{A},\rho \tilde{B})\]
where the sum is indexed by the double-coset space $\Stab A \backslash \Gamma / \Stab B$ with some restrictions analog to $\asrt{\bir>1}$ and $\asrt{\across>1}$ ensuring that it has finite support, which we shall comment on later.
These define functions on (a subset in) the space of representations $\Hom(\Gamma,\PSL_2(\Field))$ considered up to $\PSL_2(\Field)$-conjugacy at the target. One may ask for interpretations of their limiting values at special points in its appropriate compactifications.
As explained in the previous paragraph, this construction works in particular for discrete subgroups of $\PSL_2(\R)$.
In general, we may want to specify that $\rho(\Gamma)$ is a discrete subgroup of $\PSL_2(\Field)$ after $\Field$ has been given a topology, or furthermore that $\rho(\Gamma)$ has finite covolume for the Haar measure on $\PSL_2(\Field)$ with respect to a measure on $\Field$.
In that case, one may consider the quotient of the symmetric space of $\PSL_2(\Field)$ by $\rho(\Gamma)$, and observe the relative position between the "cycles" corresponding to $A,B$ in that quotient.

We may now suggest some tantalising connections between arithmetic and topology.
For this, we should compare our summations (\ref{eq:Link}) and (\ref{eq:Cos}) with the modular cocycles introduced in \cite{Duke-Imamoglu-Toth_modular-cocycles-linking_2017} and the products appearing in \cite{Darmon-Vonk_arithmetic-intersections-modular-geodesics_2022}.
Let us note however that \cite{Duke-Imamoglu-Toth_modular-cocycles-linking_2017, Matsusaka_hyperbolic-Rademacher_2023} consider the linking numbers $\lk(A+A^{-1},B+B^{-1})$ between cycles obtained by lifting a geodesic and its inverse: this number amounts to the geometric intersection $I(A,B)$ of the modular geodesics.
Furthermore \cite{Darmon-Vonk_arithmetic-intersections-modular-geodesics_2022} considers deformations of an arithmetic nature for these intersection numbers.
None of these address the actual linking numbers, and their approach is motivated by the arithmetic of modular forms, while ours will be inspired by the geometry of the character variety.
Thus it would be interesting on the one hand to understand the arithmetic of linking numbers in terms of the modular forms appearing in \cite{Katok_modular-forms-geodesics_1984} or the modular cocycles in \cite{Duke-Imamoglu-Toth_modular-cocycles-linking_2017, Matsusaka_hyperbolic-Rademacher_2023}, and on the other hand to relate the $p$-arithmetic intersections numbers considered in \cite{Darmon-Vonk_arithmetic-intersections-modular-geodesics_2022} to the special values of functions $\Link_\rho$ and $\Cos_\rho$ defined for representations $\rho \colon \PSL_2(\Z) \to \PSL_2(\Q_p)$ as suggested above.

\subsection{Special values of Poincar\'e Series}

We may apply the general averaging procedure explained in paragraph \ref{subsec:F(A,B)} to other conjugacy invariants $f_q(A,B)$ and define new functions $F_q(A,B)$ on the character variety of $\PSL_2(\Z)$.
Their limit at the boundary point $q=\infty$ will be expressed in terms of the linking number $\lk(A,B)$ as soon as $f_q(A,B)$ converges to an expression of $\cosign(A,B)$. 

Various motivations (including special values for Poincar\'e series \cite{Siegel_advanced-number-theory_1965, Dirichlet_formes-quadratiques-complexes_1842}, and McShane's identity \cite{Bowditch_McShane-Markov_1996}) suggest to choose $f_q(A,B)=(x+\sqrt{x^2-1})^{-s}$ for some variable $s\in \C$ where $x=\frac{1}{4}(\Tr(A_qB_q^{-1})-\Tr(A_qB_q))$ is the numerator of $\tfrac{1}{4}\cos(A_q,B_q)$ in the formula of Lemma \ref{Lem:cos-cosh-sinh}.
This summand $f_q(A,B)$ can also be written $e^{-si\theta}$ where $\theta$ is the angle between the oriented geometric axes of $A_q$ and $B_q$ when they intersect and $e^{-sl}$ where $l$ is the length of the ortho-geodesic arc $\gamma$ connecting the geometric axes of $A_q$ and $B_q$ when they are disjoint.
In formula:
\begin{equation*}
    F_{q}(A,B)= \sum \left(x+\sqrt{x^2-1}\right)^{-s} = \sum_{\alpha \perp \gamma \perp \beta} \exp(-sl_\gamma) - \sum_{p\in \alpha\cap \beta} \exp(-si\theta_p).
\end{equation*}
So the sum over all double cosets splits as a finite sum computable as explained in \ref{subsec:F(A,B)}, and an infinite series which converges for $\Re(s)>1$ (the topological entropy for the action of $\PSL_2(\Z)$ on the hyperbolic plane).
The infinite sum is a bivariate analog (in $(A,B)$) of the univariate Poincar\'e "theta-series" which appeared in the works of Eisenstein: those admit meromorphic continuation to $s\in \C$ and their special values in the variable $s$ have been of interest for arithmetics and dynamics. Similar Poincar\'e series associated to one modular geodesic are also defined in \cite{Katok_modular-forms-geodesics_1984}.
The earliest appearance we found for bivariate series is in \cite[Section 50]{Ford_automorphic-functions_1923}, and the only other in \cite{Paulin_series-poincare_2013}.
When $q=\infty$ and $s=1$, the real part of the finite sum evaluates to $2\lk(A,B)-I(A,B)$, but one may wonder about the infinite series (now the order in which we take limits in $s$ and $q$ may import).

More generally, one strategy to relate modular topology and quadratic arithmetic is to choose $f$ with appropriate symmetries and analyticity properties so that the sum over all double cosets can be understood: then one deduces a relationship between a topologically meaningful finite sum, and the infinite series whose special values may be of interest in arithmetic. The dilogarithm of the cross-ratio also looks like a good candidate \cite{Bridgeman_orthospectra-laminations-dilog-identities_2011}...

\subsection*{Further background references}

Besides the references already cited, let us mention a few other sources to learn some background or related notions.

The books \cite{Conway_sensual-quad-form_1997} and \cite{Hatcher_topology-numbers_2022} provide a sensual and visual introduction to the modular group, in relation with quadratic arithmetic.
The surveys \cite{Birman-Kofman_new-twist-Lorenz-links_2009, Birman_math-Lorenz-knots_2013} of what was known about Lorenz knots at that time also open new problems. 
The book \cite{Ghys-Harpe_groupes-hyperboliques-gromov_1990} and report \cite{Paulin_actions-groupes-arbres_1997} expose the theory of hyperbolic groups and their actions on trees, introduced in the seminal essay \cite{Gromov_hyperbolic-groups_1987}.
Finally the monograph \cite{Serre_arbres_1977} expounds the Bass-Serre theory of $\PSL_2$.

\bibliographystyle{alpha} 
\bibliography{paratext/bib_LiNuMoK.bib}

\end{document}